\documentclass[11pt]{article}

\usepackage[a4paper, margin=1in]{geometry}  
\usepackage{textcomp}
\usepackage{graphicx}
\usepackage{amssymb}
\usepackage{amsmath}
\usepackage[normalem]{ulem}
\usepackage{float}  
\usepackage{natbib}
\usepackage{xcolor}
\usepackage[utf8]{inputenc}
\usepackage{longtable}
\usepackage{booktabs}
\usepackage{hyperref}
\usepackage{geometry}
\usepackage{authblk}
\usepackage{amsthm}

\theoremstyle{plain}
\newtheorem{theorem}{Theorem}[section]       
\newtheorem{proposition}[theorem]{Proposition}

\newtheorem{definition}[theorem]{Definition} 

\theoremstyle{remark}
\newtheorem{remark}[theorem]{Remark}

\title{Reduced Order Data-driven Twin Models for Nonlinear PDEs by Randomized Koopman Orthogonal Decomposition and Explainable Deep Learning}
\author{D.A. Bistrian \thanks{Department of Electrical Engineering and Industrial Informatics, University Politehnica Timisoara, 300006 Timisoara, Romania; 
 diana.bistrian@upt.ro}}
\date{\today}

\begin{document}
\maketitle

\begin{abstract}
This study introduces a data-driven twin modeling framework based on modern Koopman operator theory, offering a significant advancement over classical modal decomposition by accurately capturing nonlinear dynamics with reduced complexity and no manual parameter adjustment. The method integrates a novel algorithm with Pareto front analysis to construct a compact, high-fidelity reduced-order model that balances accuracy and efficiency. An explainable NLARX deep learning framework enables real-time, adaptive calibration and prediction, while a key innovation-computing orthogonal Koopman modes via randomized orthogonal projections-ensures optimal data representation. This approach for data-driven twin modeling is fully self-consistent, avoiding heuristic choices and enhancing interpretability through integrated explainable learning techniques.
The proposed method is demonstrated on shock wave phenomena using three experiments of increasing complexity, accompanied by a qualitative analysis of the resulting data-driven twin models.
\end{abstract}


\section{Introduction}

In recent years, digital twin models have emerged as powerful tools for understanding, simulating, and predicting the behavior of complex physical systems across various domains \citep{Glaessgen2012, Grieves2017, Tao2019, RasheedDigital2020}. Despite their success in engineered systems with well-defined models and sensing infrastructures, the application of digital twins to nonlinear dynamical systems remains limited due to the inherent complexity of such dynamics. Moreover, conventional approaches often rely on full physical modeling and high computational demands, constraining their scalability and suitability for data-driven or partially observed systems.

To address key limitations in modeling nonlinear dynamical systems, this work introduces a novel theoretical formulation and an algorithmic framework within the emerging domain of data-driven twin modeling. This new paradigm is inspired by the digital twin concept but diverges fundamentally in its construction and purpose: rather than embedding a full physics-based model, the data-driven twin model is designed to learn and replicate the nonlinear behavior of a complex dynamical system directly from data, achieving high fidelity with significantly reduced model complexity. 

At the core of this contribution is a new mathematically rigorous framework that decouples system evolution into stationary and dynamic components using the principles of modal decomposition \citep{Holmes1996}.
Based on this theoretical foundation, this research presents the first mathematical introduction of a data-driven twin model (DTM) tailored with a novel algorithm, specifically formulated for unsupervised nonlinear dynamical systems. 

The proposed framework defines the DTM as a reduced-order model capable of faithfully mirroring the full system behavior with significantly reduced computational complexity. This approach addresses a critical gap in existing modeling approaches, specifically the integration of real-time adaptability with the accurate simulation of nonlinear phenomena.


\subsection{Modal Decomposition Literature Review and Previous Work}

Modal decomposition \citep{Holmes1996} represents a generic mathematical framework for the characterization of nonlinear dynamical systems in the context of reduced-order modeling  \citep{Rowley2010, Frederich2011, KutzOptical2021, Schmid2012, MezicKoop2021, Ahmed2020, XiaoFang2019}.

Proper Orthogonal Decomposition (POD) and Dynamic Mode Decomposition (DMD) are widely used modal decomposition methods for reducing computational complexity in modeling complex dynamical systems. POD, rooted in the Karhunen–Loève expansion independently developed by Loève \citep{Loeve1945} and Karhunen \citep{Karhunen1946}, was adapted to fluid mechanics by Lumley \citep{Lumley1970} to extract dominant coherent structures from turbulence. Since then, POD has been applied across turbulent and convection-dominated flows \citep{Iliescu2022, Wang2012, XiaoBis2019, Gunz2018}, engineering, and oceanography \citep{San2011, San2014, Osth2014}, often combined with data assimilation \citep{XiaoFang2018}, control \citep{SierraKaiser2021}, and machine learning, including neural networks \citep{Xiao2019b, San2019, KaptanogluBrunton2022, OwensKutz2023, ChengFang2020}. However, POD's energy-based truncation can miss low-energy but crucial modes, reducing accuracy in complex or nonlinear systems \citep{dumon2013, Dimitriu2017}, and its projection nature may fail to capture fine-scale or transient dynamics absent from training data, especially in long-term or highly nonlinear scenarios.

Schmid and Sesterhenn \citep{Schmid2008, Schmid2010} introduced Dynamic Mode Decomposition (DMD), a numerical technique based on Koopman Operator Theory \citep{Koopm1931, KoopmNeum1932}, which has been widely applied across science and engineering \citep{Rowley2010, MezicKoop2021, PercicMezic2021, Pant2021}. J. Nathan Kutz’s group significantly advanced this area by developing robust Koopman-based, data-driven modeling and system identification methods \citep{Tu2014, Brunton2016, Kutz2016, KutzOptical2021, LiKutz2022}. Many DMD variants have since emerged to enhance robustness and applicability, including optimized \citep{Chen2012}, exact \citep{Tu2014}, sparsity-promoting \citep{Jovanovic2014}, multi-resolution \citep{Kutz2016}, extended \citep{Williams2015}, recursive \citep{Noack2016}, control-enabled \citep{Proctor2016}, randomized low-rank \citep{Erichson2016}, core sketch \citep{AhmedSanBistrian2022}, bilinear \citep{Goldschmidt2021}, and higher-order DMD \citep{LeClainche2017}.

DMD captures spatiotemporal dynamics through modes with distinct frequencies and growth rates; however, its accuracy depends on proper mode selection. Identifying relevant modes often requires additional criteria involving manual tuning, domain expertise, or heuristics, increasing computational cost and limiting its use in fully unsupervised settings.

Koopman operator theory, coupled with DMD, has become a foundational tool in contemporary dynamical systems analysis and is widely regarded as a key mathematical framework of the $21^{st}$ century. The author has notably advanced Koopman-based algorithms, particularly addressing challenges posed by high-dimensional data and complex temporal dynamics.

Significant contributions include the development of mode selection strategies that leverage amplitude-weighted growth rates and energy-Strouhal number criteria to systematically identify the most dynamically relevant modes \citep{Bistrian2015, Bistrian2017a}. To overcome the computational challenges of standard DMD, the author introduced the adaptive randomized dynamic mode decomposition (ARDMD) \citep{Bistrian2017}, enhancing efficiency and scalability in reduced-order modeling, especially within fluid dynamics applications \citep{Bistrian2018}.

The techniques developed in these prior studies tackled the enduring challenge of identifying dominant modes in nonlinear and unsteady systems, which is essential for precise reconstruction and forecasting. Additionally, the author has expanded the application of DMD into emerging fields such as epidemiological modeling \citep{BistrianAIP2020, BistrianAIP2019}, demonstrating its wider societal impact.

\subsection{Problem Statement and Novelty of the Present Research}

Acknowledging the inherent limitations of classical DMD, namely, the non-orthogonality of its modes and the substantial computational overhead linked to manual mode selection, as well as the constraints of POD, which, despite yielding orthogonal modes, still depends on heuristic or user-specified criteria for mode selection, the proposed work introduces a novel approach for constructing a data-driven twin model that accurately replicates the full dynamics of nonlinear systems while significantly reducing model complexity. 

Grounded in modern Koopman operator theory, the proposed framework offers a self-consistent, interpretable, and scalable solution that eliminates manual intervention and is well-suited for unsupervised learning. By integrating explainable deep learning within a rigorous mathematical setting, it ensures transparency in the modeling process, an essential feature absent in conventional black-box machine learning models.

This novel framework significantly advances Koopman-based modal decomposition. Its key contributions include: 

$1.$ A fully self-consistent algorithm incorporating Pareto front analysis to eliminate heuristic parameter tuning;  

$2.$ Integration of explainable deep learning techniques to enable transparent, adaptive, and interpretable modeling; and

$3.$ A compact, high-fidelity data-driven twin model balancing accuracy and computational efficiency.

The structure of the article is organized as follows:
Section 2 outlines the mathematical foundations of reduced-order modeling grounded in Koopman operator theory. Section 3 addresses the computational aspects and introduces a novel modal decomposition algorithm, supported by a rigorous mathematical justification. Section 4 applies the proposed method to shock wave phenomena through three experiments of increasing complexity, while Section 5 presents the corresponding numerical results. Section 6 concludes with a summary of key findings and potential directions for future research. Appendix A provides the full derivation of the analytical and numerical solutions, and Appendix B includes the pseudocode of the proposed algorithm for reproducibility and implementation.
The key mathematical notations used in this article are summarized in Table \ref{tabnot}.
\begin{longtable}{@{}p{0.35\textwidth}p{0.6\textwidth}@{}}
\caption{Key notations used throughout the article.\label{tabnot}} \\
\toprule
\textbf{Mathematical notation} & \textbf{Description} \\
\midrule
\endfirsthead

\toprule
\textbf{Mathematical notation} & \textbf{Description} \\
\midrule
\endhead

\midrule
\multicolumn{2}{r}{\small\itshape Continued on next page} \\
\endfoot

\bottomrule
\endlastfoot

$\mathcal{K}^t$ & Koopman operator at time $t$, indicating its dependence on the temporal variable in continuous-time dynamical systems. \\
$\{\mathcal{K}^t\}_{t \geq 0}$ & The semigroup of Koopman operators. \\
$\mathcal{K}$ & Koopman operator for discrete-time dynamical system. \\
$N_x$ & Number of spatial grid points. \\
$N_t$ & Number of time steps or snapshots, excluding the initial state. \\
$\left\| \cdot \right\|$ & The norm for functions in the Hilbert space $L^2(\Omega)$. \\
$\left\| \cdot \right\|_2$ & The Euclidean (or $\ell^2$) norm for vectors. \\
$u(x,t)$ & Solution of the PDE, representing the system state at spatial location $x$ and time $t$. \\
$u_i = u(x,t_i),\,i = 1, \dots, N_t$ & Measurements of the PDE solution sampled uniformly in time with step size $\Delta t$. \\
$u_i^{(k)} = u^{(k)}(x,t_i)$ & Reduced-order approximation of the PDE solution at time $t_i$, using the $k$-th order model. \\
$N_{DTM}$ & Rank of the optimal Koopman basis, computed by the algorithm. \\
$\phi_j(x),\, j = 1, \ldots, N_{DTM}$ & Leading Koopman modes computed by the algorithm, representing dominant spatial patterns. \\
${a_j}\left( {{t_i}} \right) $ & The modal growing amplitudes.\\
$\hat a(t)$ & The system temporal dynamics in the reduced-order space, modeled by the deep learning architecture. \\
$u^{DTM}(x,t)$ & The reduced-order data-driven twin model for the PDE solution. \\
$H$ & Conjugate transpose (Hermitian transpose); for a matrix $A$, $A^H = \overline{A}^T$. \\
$\langle \cdot, \cdot \rangle$ & Inner product operator; for vectors $x, y \in \mathbb{C}^n$, $\langle x, y \rangle = x^H y$. \\
$ X_{.,j} $ & The \(j\)th column of the matrix $X$.\\
\end{longtable}

\section{Mathematical Background}

\subsection{Computational Complexity Reduction of PDEs by Reduced-order Modeling}

Partial differential equations (PDEs) serve as a fundamental framework for accurately modeling nonlinear dynamical phenomena, owing to their capacity to represent spatial-temporal variations in complex systems. Nevertheless, the numerical solution of PDEs is often computationally demanding, particularly in high-dimensional and strongly nonlinear settings.
This computational burden necessitates strategies for complexity reduction to enable practical simulation and analysis. 
In this context, by reducing the dimensionality of the system while retaining its essential dynamics, reduced-order modeling significantly lowers computational cost and facilitates efficient simulation and analysis.

Suppose that $\Omega \subset \mathbb{R}^d$ represents the bounded computational domain and let the Hilbert space ${L^2}\left( \Omega \right)$  of square integrable functions on $\Omega$:
\begin{equation}
{L^2}\left( \Omega  \right) = \left\{ {\phi :\Omega  \to \mathbb{R}\left| {\int_\Omega  {{{\left| \phi  \right|}^2}d\Omega  < \infty } } \right.} \right\}
\end{equation}
to be endowed with the inner product:
\begin{equation}
\left\langle {{\phi _i},{\phi _j}} \right\rangle  = \int_\Omega  {{\phi _i}} \overline {{\phi _j}} \;d\Omega \quad for\;{\phi _i},{\phi _j} \in {L^2}\left( \Omega  \right)
\end{equation}
and the induced norm $\left\| \phi  \right\| = \sqrt {\left\langle {\phi ,\phi } \right\rangle }  $  for $\phi  \in {L^2}\left( \Omega \right)$.

Consider a nonlinear PDE for state variable $u\left( {x,t} \right)$, $x \in {\Omega _x} \subset {\mathbb{R}^d}$, $t \in {\Omega _t} = \left[ {{t_0},T} \right]$,
on domain $D = {\Omega _x} \times {\Omega _t}$, of the form:
\begin{equation}\label{pde}
    \frac{\partial u(x,t)}{\partial t} = \mathcal{N}[u(x,t)],
\end{equation}
 where \( \mathcal{N} \) is a nonlinear differential operator, 
with the initial condition: 
\begin{equation}
IC:\quad u\left( {x,{t_0}} \right) = {u_0}\left( x \right).
\end{equation}
The Robin boundary conditions, defined by two time-continuous functions ${g_a}\left( t \right)$ and ${g_b}\left( t \right)$, specified at the endpoints of the spatial domain, take the general form:
\begin{equation}
BC:\quad \left\{ {\begin{array}{*{20}{l}}
{{\alpha _1}u\left( {a,t} \right) + {\alpha _2}{u_x}\left( {a,t} \right) = {g_a}\left( t \right)},\\
{{\beta _1}u\left( {b,t} \right) + {\beta _2}{u_x}\left( {b,t} \right) = {g_b}\left( t \right)},
\end{array}} \right.
\end{equation}
where ${\alpha _1}$, ${\alpha _2}$, ${\beta _1}$, ${\beta _2}$ are real prescribed coefficients.

Over the past half-century, partial differential equations have been predominantly solved using numerical techniques such as the finite difference and finite element methods. While these computational tools have significantly advanced our ability to approximate complex PDE solutions, they also introduce substantial implementation complexity and often require a solid understanding of the underlying analytical structure of the equations. In modern computational mathematics, reduced-order modeling has emerged as a powerful approach to mitigate computational demands by simplifying high-dimensional PDE systems, enabling more efficient and tractable solutions without sacrificing essential accuracy.

\begin{definition}
The principle of reduced-order modeling aims to find an approximate solution of the form
\begin{equation}\label{rom1}
\left\{ {\begin{array}{*{20}{l}}
{u\left( {x,t} \right) = \mathop {\lim }\limits_{p \to \infty } \sum\limits_{j = 1}^p {{a_j}\left( t \right){\psi _j}\left( x \right)} ,\quad t \in {\mathbb{R}_{ \ge 0}}}\\
{\frac{{d{a_j}\left( t \right)}}{{dt}} = \mathcal{M}\left( {{a_j}\left( t \right),t} \right),\;{a_j}\left( {{t_0}} \right) = a_j^0},
\end{array}} \right.
\end{equation}
where \( \{ \psi_j(x) \}_{j=1}^{p} \subset L^2(\Omega) \) is a set of basis functions and \( \mathcal{M} \) is a locally Lipschitz continuous map governing the temporal dynamics of the coefficients \( a_j(t) \). The initial condition is prescribed on the coefficient vector at the initial time \(t_0\), with \(a_j(t_0) = a_j^0\).

The series is assumed to converge in the \( L^2(\Omega) \) norm for each fixed time \( t \), that is,
\[
\lim_{p \to \infty} \left\| u(x,t) - \sum_{j=1}^p a_j(t) \psi_j(x) \right\|= 0.
\]

This reduced-order model aims to approximate the true solution of the full-order system defined by the PDE (\ref{pde}), while ensuring a small global approximation error and computational stability.

Expression (\ref{rom1}) is called the \emph{reduced-order model} of the PDE (\ref{pde}).
\end{definition}

The basis functions $\Psi$ are called modes or shape-functions, therefore the representation (\ref{rom1}) is called modal decomposition.
Two primary goals are achieved through this modal decomposition approach: (1) the decoupling of spatial structures from the unsteady (time-dependent) components of the solution, and (2) the modeling of the original system’s nonlinearity through a lower-dimensional temporal dynamical system. This separation not only reduces computational complexity but also preserves the essential dynamics of the original PDE system in a more efficient and interpretable form.

\subsection{Koopman Operator Framework for Modal Decomposition}

Introduced by Koopman in 1931, Koopman Operator Theory \citep{Koopm1931, KoopmNeum1932, Koopm1936} provides a rigorous mathematical framework for the modal decomposition of nonlinear dynamical systems by representing them through linear, infinite-dimensional operators acting on observables. Despite its early formulation, the theory remained largely theoretical for decades.

In 1985, Lasota and Mackey \citep{LasotaMackey1985} brought renewed scientific research interest to the original approach proposed by Koopman. Notably, they were the first to explicitly introduce and formalize the term "Koopman operator" in the literature, thereby establishing a standardized nomenclature for this concept. This terminology was further reinforced in the second edition of their book \citep{LasotaMackey94}, where the use of the Koopman operator as the adjoint to the Frobenius-Perron operator was maintained and elaborated upon functional analysis in space ${L^\infty }$.

However, at that time, computational resources were still insufficient for the numerical treatment of large-scale complex systems.
The practical potential of Koopman Operator Theory began to be realized in the early 2000s.

Its modern scientific revival began with the pioneering work of Mezi\'{c} and Banaszuk \citep{MezicBanaszuk2004, Mezic2005} in fluid dynamics applications, establishing the foundation for practical, data-driven modal decomposition analysis.
The seminal work by Mezi'{c} \citep{Mezic2005} introduced a rigorous spectral framework for analyzing nonlinear dynamical systems via Koopman operator theory, formalizing the spectral decomposition of observables into Koopman eigenfunctions and modes, and establishing the foundation for operator-theoretic analysis of nonlinear dynamics. The spectral properties of the Koopman operator have been the subject of ongoing investigation, as exemplified in \citep{Rowley2009, Chen2012, Brunton2016}.

In this paper, we propose an adaptation of the Koopman framework to the context of nonlinear partial differential equations, extending its applicability to infinite-dimensional dynamical systems governed by PDEs.

Let \( \mathcal{U}(t) \) denote the evolution operator of (\ref{pde}) such that:
\begin{equation}
    u(x,t) = \mathcal{U}(t) u(x,{t_0}),
\end{equation}
i.e., \( \mathcal{U}(t) \) maps the initial condition \( u(x,{t_0}) \) to the state at time \( t \).

\begin{definition}
An observable of the PDE (\ref{pde}) is defined as a possibly nonlinear functional on the state space, typically a Hilbert space such as  $\mathcal{H}={L^2}\left( \Omega  \right)$:
\begin{equation}
 g: \mathcal{H} \to \mathbb{C}^m,
\end{equation}
which maps the infinite-dimensional state $u$ to a measurable output quantity. Here $m\in \mathbb{N}$ denotes the dimension of the observable space, representing the number of scalar measurements or features extracted from the state $u$. These measurements may vary over time, thus capturing the temporal dynamics of the system in a finite-dimensional representation.
\end{definition}

\begin{remark}\label{rem1}
The choice of observable \( g \) in the Koopman framework is flexible and can range from simple to highly structured representations. A common example is the \emph{identity observable}, where \( g(u) = u \), mapping the system state directly to itself. Alternatively, one may consider an infinite-dimensional observable vector, comprising a set of functions of the state, such as
\begin{equation}
g(u) = {[u,\;\kappa u,\;{u^2},\;{u^3},\;\sin \left( u \right), \ldots ]^T},
\end{equation}
where \( \kappa \) denotes a linear operator (e.g., a differential operator), and the higher-order terms capture nonlinear features of the system. This extended observable space facilitates a richer linear representation of the underlying nonlinear dynamics through the Koopman operator.
\end{remark}

\begin{definition}
The semigroup of Koopman operators \( \{\mathcal{K}^t \}_{t \geq 0}\) acts on observable functions by composition with the evolution operator of the states:
\begin{equation}
    \left( \mathcal{K}^t g \right)({u_0}\left( x \right)) = g\left( \mathcal{U}(t) {u_0}\left( x \right) \right),
\end{equation}
where \( {u_0}\left( x \right) = u\left( {x,{t_0}} \right) \) is the initial condition.
As a result, the Koopman operator is also known as a composition operator \citep{Koopm1931}.
\end{definition}

In practice, the Koopman operator is considered to be the dual, or more precisely, the left adjoint, of the Ruelle-Perron-Frobenius transfer operator \citep{Perron1907, Frobenius1908, Ruelle1993}, which acts on the space of probability density functions.

\begin{remark}
The Koopman operator is \textit{linear}, even though the underlying dynamics may be nonlinear \citep{Koopm1931}:
\begin{equation}
    \mathcal{K}^t (a g_1 + b g_2) = a \mathcal{K}^t g_1 + b \mathcal{K}^t g_2.
\end{equation}
\end{remark}

\begin{remark}
The infinitesimal generator \( \mathcal{L} \) of the Koopman semigroup governs the evolution of observables via
\begin{equation}
    \frac{d}{dt} g(u) = \mathcal{L} g(u),
\end{equation}
where, for deterministic dynamics $\dot u = \mathcal{N}\left[ u \right]$, the generator acts as
\begin{equation}
    \mathcal{L} g(u) = \left\langle \frac{\delta g}{\delta u},\, \mathcal{N}[u] \right\rangle,
\end{equation}
with \( \frac{\delta g}{\delta u} \) denoting the functional (Fréchet) derivative of \( g \) with respect to \( u \). The inner product is taken over the spatial domain \citep{LasotaMackey1985}.
\end{remark}

Following Koopman operator theory, the solution \( u(x,t) \) can be expressed as a modal decomposition of the form:
\begin{equation}
    u(x,t) = \sum_{j=1}^\infty a_j(t) \, \Phi_j(x),
\end{equation}
where \( \Phi_j(x) \) are spatial Koopman modes, and \( a_j(t) \) are time-dependent coefficients.

The evolution of the modal coefficients \( \mathbf{a}(t) = [a_1(t), a_2(t), \dots]^T \) is governed by a (possibly nonlinear) dynamical system:
\begin{equation}
    \frac{d \mathbf{a}(t)}{dt} = \mathcal{M}(\mathbf{a}(t), t),
\end{equation}
where \( \mathcal{M} \) is a nonlinear operator (or vector field) that defines the temporal evolution of the coefficients \( \mathbf{a}(t) \).
Thus, the reduced-order solution of form (\ref{rom1}) to the PDE is achievable form the perspective of Koopman mode decomposition.

Koopman modes are generally not unique, and their characterization depends on several interrelated factors. First, the definition of Koopman modes is inherently linked to the choice of observables, as the Koopman operator acts on a space of functions rather than directly on the state space. Different selections of observables can lead to different spectral representations, thereby affecting the resulting modes. 

This intrinsic non-uniqueness is further compounded in practical settings where Koopman modes are often approximated using numerical techniques such as DMD \citep{Schmid2008, Schmid2010}. Such methods introduce additional variability due to factors including data noise, finite sample sizes, and specific algorithmic implementations. Moreover, in cases where the Koopman operator exhibits a degenerate spectrum \citep{MezicKoop2021, Mezic2022}, the associated eigenspaces are typically multidimensional, allowing for multiple valid decompositions and further contributing to the non-uniqueness of the Koopman modes. 

This work proposes a new algorithm, rooted in Koopman operator theory but fundamentally different in its formulation, specifically designed for the construction of a data-driven twin model. This new method computes Koopman modes and their associated temporal coefficients by leveraging the Koopman propagator operator projected onto a Krylov subspace of optimally selected rank. In contrast to classical DMD, the proposed approach facilitates a more structured and computationally efficient extraction of spectral components by constraining the dynamics within a finite-dimensional subspace tailored to capture the most relevant features of the system's evolution.

\section{Computational Aspects }

This section introduces a novel methodology for constructing a data-driven twin model of a partial differential equation using Koopman mode decomposition. The proposed approach comprises two primary stages. The first stage, referred to as the offline phase, involves the extraction of Koopman modes and the corresponding modal decomposition coefficients that characterize the state evolution of the system. The second stage, known as the online phase, focuses on identifying a nonlinear continuous-time dynamical model that governs the temporal behavior of the data-twin, thereby enabling real-time simulation and prediction. A detailed account of the computational implementation of the proposed methodology is presented hereinafter.

\subsection{Offline Phase:  Koopman Randomized Orthogonal Decomposition Algorithm}

The data ${u_i}= u\left( {x,{t_i}} \right),\;{t_i} = i\Delta t,\;i = 0,...,{N_t}$, represent measurements of the PDE solution at the constant sampling time $\Delta t$, $x$ representing the Cartesian spatial coordinate.

The data matrix whose columns represent the individual data samples is called \textit{the snapshot matrix}
\begin{equation}\label{snap}
V = \left[ {\begin{array}{*{20}{c}}
{{u_0}}&{{u_1}}&{...}&{{u_{N_t}}}
\end{array}} \right] \in {\mathbb{R}^{{N_x} \times ({N_t} + 1)}}.
\end{equation}
Each column ${u_i} \in {\mathbb{R}^{N_x}}$ is a vector with ${N_x}$ components, representing the spatial measurements corresponding to the ${N_t} + 1$ time instances.

We aim to find the data-driven twin  representation of the solution at every time step $\left\{{{{\rm{t}}_1}{\rm{,}}...{\rm{,}}{{\rm{t}}_{N_t}}} \right\}$,  according to the following relation:
\begin{equation}\label{dtm}
{u_i} = u\left( {x,{t_i}} \right) = \sum\limits_{j = 1}^{{N_{DTM}}} {\underbrace {{a_j}\left( {{t_i}} \right)}_{Modal\;amplitudes}\overbrace {{\phi _j}\left( x \right)}^{Leading\;Koopman\;modes},} \;\;i = 1,...,{N_t},
\end{equation}
where $\Phi  = \left\{ {{\phi _j}} \right\}_{j = 1}^{{N_{DTM}}}$ represent the extracted Koopman modes base functions, which we call \textit{the leading Koopman modes}, ${N_{DTM}} \ll \min \left( {{N_x},{N_t} + 1} \right)$ represents the number of terms in the representation (\ref{dtm}) which we impose to be minimal and ${a_j}\left( {{t_i}} \right) $ represent the modal growing amplitudes.

In classical Dynamic Mode Decomposition the Koopman modes are not orthogonal, which often necessitates a large number of modes to achieve an accurate modal decomposition. This lack of orthogonality represents a significant drawback, as it may lead to redundancy and reduced computational efficiency. 

It is more convenient to seek an orthonormal base of Koopman modes
\begin{equation}
 \Phi  = \left\{ {{\phi _1},{\phi _2},...} \right\},\quad \left\langle {{\phi _i}\left( x \right),{\phi _j}\left( x \right)} \right\rangle  = {\delta _{ij}},\quad \left\| \phi  \right\| = 1,
\end{equation}
where ${\delta _{ij}}$ is the Kronecker delta symbol, consisting of a minimum number of Koopman modes ${{\phi _i}\left( x \right)}$, such that the approximation of $u\left( {x,t} \right)$ through this base is as good as possible, in order to create a data-twin model to the PDE, of reduced computational complexity.

The objective of the present technique is to represent with highest precision the original data through the data-twin model (\ref{dtm}) having as few terms as possible. To achieve this goal, several major innovative features are introduced in the proposed algorithm:
\begin{itemize}
\item	 The method yields orthonormal Koopman modes, guaranteeing mutual orthogonality and thereby enabling a more compact representation of the system dynamics (see Theorem~\ref{th1} and its proof).
\item	To avoid the computational burden associated with traditional high-dimensional algorithms, the proposed method employs a randomized singular value decomposition (RSVD) technique for dimensionality reduction. The incorporation of randomization offers a key advantage: it eliminates the need for additional selection criteria to identify shape modes, which is typically required in classical approaches such as DMD or POD. The resulting algorithm efficiently identifies the optimal reduced-order subspace that captures the dominant Koopman mode basis, thereby ensuring both computational efficiency and representational fidelity (see Appendix~\ref{appB}).
\item The methodology aims to achieve maximal correlation and minimal reconstruction error between the data-driven twin model and the exact solution of the governing PDE (see Section~\ref{numres}).
\end{itemize}

The theoretical foundation of RSVD was laid in the seminal work of Halko, Martinsson, and Tropp (2011) \citep{Halko2011} as part of a broader class of probabilistic algorithms for matrix approximation. Subsequent study of Erichson and Donovan (2016), demonstrated its practical efficiency in motion detection \citep{Erichson2016}. Bistrian and Navon (2017) extended its use for the first time in fluid dynamics \citep{Bistrian2017, Bistrian2018}.
Recently, the theory of randomization has been further employed in the development of reduced-order modeling techniques, particularly in conjunction with projection learning methods \citep{Surasinghe2021}, block Krylov iteration schemes \citep{Qiu2024} or deep learning \citep{Bistrian2022}.

\begin{proposition} \label{Randsvd} ($k$-RSVD: Randomized Singular Value Decomposition of rank $k$)

Let $V_0 \in \mathbb{R}^{N_x \times N_t}$ be a real-valued data matrix, and let $k \in \mathbb{N}$ with $k < \min(N_x, N_t)$ be the target rank.

Then, the rank-$k$ Randomized Singular Value Decomposition ($k$-RSVD) of $V_0$ yields an approximate factorization:
\[
V_0 \approx U \Sigma W^H,
\]
where:
$U \in \mathbb{R}^{N_x \times k}$ contains approximate left singular vectors,

$\Sigma \in \mathbb{R}^{k \times k}$ is diagonal with approximate singular values,

$W \in \mathbb{R}^{N_t \times k}$ contains approximate right singular vectors.

This decomposition is computed via the $k$-RSVD algorithm described in Appendix \ref{appB}.
\end{proposition}

\begin{theorem} (Koopman Randomized Orthogonal Decomposition)\label{th1}

Let $\{ u_i \}_{i=0}^{N_t} \subset \mathbb{R}^{N_x}$ denote a sequence of solution snapshots of a PDE, arranged into two time-shifted data matrices defined as:
\begin{equation}
\begin{array}{l}
{V_0} = \left[ {\begin{array}{*{20}{c}}
{{u_0}}&{{u_1}}&{...}&{{u_{{N_t} - 1}}}
\end{array}} \right]  \in {\mathbb{R}^{{N_x} \times {N_t}}},\;\\
{V_1} = \left[ {\begin{array}{*{20}{c}}
{{u_1}}&{{u_2}}&{...}&{{u_{{N_t}}}}
\end{array}} \right]  \in {\mathbb{R}^{{N_x} \times {N_t}}}.
\end{array}
\end{equation}

Let $k \in \mathbb{N}$ with $k < N_t$ be a target rank. Then there exists a reduced $k$-order approximation $u_i^{(k)}$ of the full state $u_i$ such that:
\begin{equation}\label{RKOD}
u_i^{(k)} = u^{(k)}(x, t_i) = \sum_{j=1}^{k} a_j(t_i) \, \phi_j(x), \quad \forall i \in \{1, \dots, N_t\},
\end{equation}
where $\{\phi_j\}_{j=1}^k$ is an orthonormal set of Koopman modes spanning the reduced subspace:
\begin{equation}\label{modesdef}
\mathrm{span}\left\{ \phi_1, \phi_2, \dots, \phi_k \right\}, \quad 
\phi_j = \frac{\langle U, X_{.,j} \rangle}{\left\| \langle U, X_{.,j} \rangle \right\|_2}, \quad j = 1, \dots, k.
\end{equation}
Here:

$U \in \mathbb{R}^{N_x \times k}$ is the matrix of the first $k$ left singular vectors of $V_0$ obtained via $k$-RSVD,

$X \in \mathbb{R}^{k \times k}$ contains the eigenvectors of the Gram matrix associated with the finite-dimensional approximation of the Koopman operator $\mathcal{K}$ defined by the time evolution:
    \[
    u_{N_t} = \mathcal{K}^{N_t} u_0.
    \]

The corresponding modal coefficients $a_j(t_i)$ are computed via the inner product between the data matrix $V_0$ and Koopman modes:
\begin{equation}\label{ampldef}
a_j(t_i) = \langle V_0, \phi_j \rangle, \quad \forall j \in \{1, \dots, k\}, \quad \forall i \in \{1, \dots, N_t\}.
\end{equation}
\end{theorem}

The proposed modal representation method is hereafter referred to as Koopman Randomized Orthogonal Decomposition (KROD).

\begin{proof}

We consider a sequence of solution snapshots $\{ u_i \}_{i=0}^{N_t} \subset \mathbb{R}^{N_x}$ obtained from a discrete-time dynamical system governed by an underlying (possibly nonlinear) evolution law. Under the Koopman operator framework \citep{Koopm1931}, we adopt the \emph{identity observable} $g(u) = u$, so that the state vectors themselves are treated as observables (see Remark~\ref{rem1}).

There exists a linear Koopman operator $\mathcal{K}$ acting on observables such that:
\begin{equation}\label{steps1}
u_i = \mathcal{K}^{i} u_0, \quad \text{for } i = 0, 1, \dots, N_t.
\end{equation}

This implies that each state $u_i$ is obtained by iterated application of the Koopman operator to the initial state $u_0$ under the identity observable. The operator $\mathcal{K}$ thus advances the system in time linearly in the (possibly infinite-dimensional) observable space.

The sequence (\ref{steps1}) lies in the Krylov subspace 
\begin{equation}\label{Krylov}
{\mathcal K_{{N_t} + 1}}\left( {\mathcal K,{u_0}} \right) = span\left\{ {{u_0},\mathcal K{u_0},{\mathcal K^2}{u_0},...,{\mathcal K^{{N_t}}}{u_0}} \right\},
\end{equation}
that captures how the observable ${u_0}$ evolves under the Koopman dynamics.

The Krylov subspace (\ref{Krylov}) gives a finite-dimensional proxy to study the action of the infinite-dimensional Koopman operator.
For a sufficiently long sequence of the snapshots, suppose that the last snapshot $u_{N_t}$ can be written as a linear combination of previous
${N_t}$ vectors, on the Krylov subspace (\ref{Krylov}),  such that:
\begin{equation}\label{linearcomb}
u_{N_t} = {c_0}{u_0} + {c_1}{u_1} + ... + {c_{{N_t} - 1}}{u_{{N_t} - 1}} + \mathcal{R},
\end{equation}
in which ${{\rm{c}}_i} \in {\rm{\mathbb{R},i  =  0,}}...{\rm{,{N_t}  -  1}}$  and $\mathcal{R}$ is the residual vector.  The following relations are true:
\begin{equation}\label{kryl}
\left\{ {{u_1},{u_2},...{u_{N_t}}} \right\} = {\mathcal K}\left\{ {{u_0},{u_1},...{u_{{N_t} - 1}}} \right\} = \left\{ {{u_1},{u_2},...,V_0c} \right\} + {\mathcal R},
\end{equation}
where $c = {\left( {\begin{array}{*{20}{c}} {{c_0}}&{{c_1}}&{...}&{{c_{{N_t} - 1}}}\end{array}} \right)^T}$ is the unknown column vector.

Assume there exists a diagonal matrix $\mathcal{S}$ whose eigenvalues asymptotically approximate those of $\mathcal{K}$ as ${\left\| {\mathcal R} \right\|_2} \to 0$ approaches zero.
Solving the eigenvalue problem:
\begin{equation}\label{approx}
V_1 = \mathcal{K}V_0 = V_0\mathcal{S} + \mathcal{R},
\end{equation}
is equivalent to solve the minimization problem:
\begin{equation}\label{minpro}
\mathop {min}\limits_\mathcal{S} \;\mathcal{R} = {\left\| {V_1 - V_0\mathcal{S}} \right\|_2},
\end{equation}
where ${\left\| {\, \cdot \,} \right\|_2}$ is the ${L_2}$-norm of ${\mathbb{R}^{{N_x} \times {N_t}}}$.

For $1 \le k \le {\rm{min}}\left( {{N_x},{N_t}} \right)$,
we identify the $k$-RSVD of $V_{0}$, that yields the factorization:
\begin{equation}
V_{0}= U\Sigma {W^H},
\end{equation}
where $U \in {\mathbb{R}^{{N_x} \times {k}}}$ and $W  \in {\mathbb{R}^{{N_t} \times {k}}}$ are orthogonal matrices that
contain the eigenvectors of ${{\rm{V}}_0}{{\rm{V}}_0}^H$ and ${{\rm{V}}_0}^H{{\rm{V}}_0}$, respectively,
 $\Sigma  = diag\left( {{\sigma _1},...,{\sigma _{k}}} \right) \in {\mathbb{R}^{k \times k}}$
is a square diagonal matrix containing the singular values of $V_{0}$ and ${H}$ means the conjugate transpose.

Relations $\mathcal{K}V_0= V_1= V_0S + \mathcal{R}, {\left\| \mathcal{R} \right\|_2}
\to 0$ and $V_0 = U\Sigma {W^H}$ yield:
\[\mathcal{K}U\Sigma {W^H} = V_1 = U\Sigma {W^H}\mathcal{S}\; \Rightarrow \; {U^H}\mathcal{K}U\Sigma {W^H} = {U^H}U\Sigma {W^H}\mathcal{S}\; \Rightarrow  \;\mathcal{S} = {U^H}\mathcal{K}U.     \]

As a consequence, the solution to the minimization problem (\ref{minpro}) is the matrix operator:
\begin{equation}\label{defS}
  \mathcal{S} = {U^H}\left( {V_1W{\Sigma ^{ - 1}}} \right).
\end{equation}

As a direct result, the eigenvalues of  $\mathcal{S}$ will converge toward the eigenvalues of the Koopman propagator operator $\mathcal{K}$.

Let $X \in {\mathbb{C}^{{k} \times {k}}}$,
$\Lambda  \in {\mathbb{C}^{{k} \times {k}}}$
be the eigenvectors, respectively the eigenvalues of the Hermitian positive semi-definite Gram  matrix $\mathcal{G} = {\mathcal{S}^H}\mathcal{S}$:
\begin{equation}\label{Geig}
  \mathcal{G}X=X\Lambda.
\end{equation}

From the $k$-RSVD of $V_0$, it follows that the column space of $V_0$ is the same as the column space of $U$, i.e.,
\begin{equation}\label{range1}
\text{Range}(V_0) = \text{Range}(U).
\end{equation}

Using the definition $V_1 = V_0 \mathcal{S}$ and substituting into the expression for $\mathcal{S}$:
\begin{align*}
\mathcal{S} &= U^H V_1 W \Sigma^{-1} \\
  &= U^H (V_0 \mathcal{S}) W \Sigma^{-1} \\
  &= U^H (U \Sigma W^H \mathcal{S}) W \Sigma^{-1} \\
  &= (U^H U) \Sigma (W^H \mathcal{S} W) \Sigma^{-1} \\
  &= \Sigma (W^H \mathcal{S} W) \Sigma^{-1},
\end{align*}
since $U^H U = I$.

Since $\mathcal{S}$ is diagonal and $W$ is unitary, $W^H \mathcal{S} W$ represents a unitary similarity transform \citep{Golub2013} of matrix  $\mathcal{S}$, i.e. $\mathcal{S}$ acts in a reduced coordinate system defined by the $k$-RSVD.

Let
\begin{equation}\label{defphi}
\Phi = U X.
\end{equation}
Since the columns of $X$ are vectors in ${\mathbb{C}^{{k} }}$ and $U$ is an orthonormal basis for the column space of $V_0$, the columns of $\Phi$ are linear combinations of the columns of $U$.

Thus,
\begin{equation}\label{range2}
\text{Range}(\Phi) \subseteq \text{Range}(U) = \text{Range}(V_0),
\end{equation}
i.e. the set of vectors $\Phi = U X$ lies in the column space of $V_0$. Therefore, the span of the eigenvectors $\Phi$ forms a subspace of the space containing the data in $V_0$.

It follows that there exists a subspace of the space containing the data in $V_{0}$, spanned by the sequence of Koopman modes:
\begin{equation}\label{modesdef}
span\left\{ {{\phi _1},{\phi _2},...,{\phi _{k}}} \right\}, \; \left\{ {{\phi _j}} \right\}_{j = 1}^k = \left\langle {U,{X_{.,j}}} \right\rangle,
\end{equation}
where $U$ represents the matrix of left singular vectors produced by $k$-RSVD of data matrix $V_{0}$. 

We proceed to demonstrate that the columns of  \(\Phi = U X\) are orthonormal.
Since \(U\) comes from the $k$-RSVD of \(V_0\), it has orthonormal columns:
\begin{equation}
U^H U = I.
\end{equation}
Now compute the inner product of \(\Phi\):
\begin{equation}
{\Phi ^H}\Phi  = {(UX)^H}(UX) = {X^H}{U^H}UX = {X^H}X.
\end{equation}

So the orthonormality of \(\Phi\) depends on \(X^H X\). The matrix \(X\) contains the eigenvectors of the Hermitian Gram matrix \(\mathcal{G} = \mathcal{S}^H \mathcal{S}\). Since \(\mathcal{G}\) is Hermitian, it has a complete set of orthonormal eigenvectors. Thus, by choosing the eigenvectors to be orthonormal, we have:
\begin{equation}
{X^H}X = I\quad  \Rightarrow \quad {\Phi ^H}\Phi  = {X^H}X = I,
\end{equation}
i.e. columns of $\Phi$ are orthonormal.

It follows that:
\begin{equation}
{\left\langle {{\phi _i},{\phi _j}} \right\rangle }= {\delta _{ij}},\quad 1 \le i \le j \le {k}
\end{equation}
i.e., Koopman modes $\Phi$ form an orthnormal base to the data space.

We continue by establishing the formula for the modal coefficients.

Assume that the data evolves linearly in reduced $k$-order representation (\ref{RKOD}), written in matrix formulation:
\begin{equation}
V_1 = \Phi A.
\end{equation}
Multiplying both sides on the left by \( \Phi^H \) results:
\begin{equation}
\Phi^H V_1 = \Phi^H \Phi A.
\end{equation}
Since \( \Phi \) has orthonormal columns, we get:
\begin{equation}
\Phi^H \Phi = I
\quad \Rightarrow \quad \Phi^H V_1 = A.
\end{equation}
Thus:
\begin{equation}
A = \Phi^H V_1.
\end{equation}
If \( V_1 = V_0 \mathcal{S} \), and \( \Phi = U X \), then we can alternatively write:
\begin{equation}
A = {\Phi ^H}{V_1} = {\Phi ^H}{V_0}\mathcal{S}.
\end{equation}

We choose to project \( V_0 \) onto the subspace spanned by the orthonormal basis \( \Phi \) to reduce computational complexity. Consequently, the coefficient matrix
\begin{equation}
A = \Phi^H V_0
\end{equation}
provides a representation with a smaller approximation error in the least-squares sense.
\end{proof}

\begin{proposition}(Multi-objective twin-model selection via Pareto front analysis) \label{optim_model_selection}

Let $\{k_1, k_2, \dots, k_N\} \subset \mathbb{N}$ be a strictly increasing sequence of target ranks satisfying
\[
1 < k_1 < k_2 < \dots < k_N \leq N_t.
\]

For each target rank $k_i$, define a Koopman-based model triplet $T_i = (\Phi_i, A_i, k_i)$ consisting of:

$\Phi_i$: a Koopman mode basis of size $k_i$,

$A_i$: the corresponding modal amplitude matrix,

$k_i$: the reduced-order rank.

Define the candidate set of model triplets:
\[
\mathcal{T} = \{ T_1, T_2, \dots, T_N \},
\]
and let $u^{(T_i)}$ denote the reduced-order solution associated with $T_i$, given by:
\[
u^{(T_i)}(x, t) = \sum_{j=1}^{k_i} a_j(t)\, \phi_j(x), \quad \text{for all } t \in \{t_1, \dots, t_{N_t} \}.
\]

Given the full-order PDE solution $u$, the goal is to identify the optimal triplet by solving the following constrained multi-objective optimization problem:
\begin{equation}\label{optimprob}
\begin{aligned}
& \text{Find } \mathcal{T}^* \subseteq \mathcal{T}, \quad \text{such that:} \\
& \text{Objective 1:} \quad \min_{T_i \in \mathcal{T}} \; E(T_i) := \| u - u^{(T_i)} \|_2, \\
& \text{Objective 2:} \quad \max_{T_i \in \mathcal{T}} \; C(T_i) := \frac{ \| u^H u^{(T_i)} \|_2^2 }{ \| u^H u \|_2 \cdot \| (u^{(T_i)})^H u^{(T_i)} \|_2 },
\end{aligned}
\end{equation}
where $\| \cdot \|_2$ denotes the discrete Euclidean norm over all space-time points, and $H$ denotes the Hermitian (conjugate transpose).

The Pareto-optimal subset of model triplets is defined as:
\begin{align*}
\mathcal{T}^* = \big\{ T_i \in \mathcal{T} \;\big|\; \nexists T_j \in \mathcal{T} \text{ such that } &\; E(T_j) \le E(T_i), \\
& C(T_j) \ge C(T_i),\; \text{and at least one inequality is strict} \big\}.
\end{align*}

Let the Pareto-optimal set of Koopman model triplet be denoted by the  solution that represent the best trade-offs between competing objectives:
\[
\mathcal{T}^* = (\Phi_{\mathrm{DTM}}, A_{\mathrm{DTM}}, N_{\mathrm{DTM}}),
\]
where the number of modes in the final reduced model is $N_{\mathrm{DTM}}$.

The resulting data-driven twin model (DTM) is then expressed as:
\begin{equation}\label{phase1dtm}
u_i^{\mathrm{DTM}}(x) = \sum_{j=1}^{N_{\mathrm{DTM}}} a_j(t_i)\, \phi_j(x), \quad i = 1, \dots, N_t,
\end{equation}
where $a_j(t_i) = (A_{\mathrm{DTM}})_{i,j}$ and $\phi_j(x) = (\Phi_{\mathrm{DTM}})_{:,j}$.
\end{proposition}

This framework is designed to simultaneously balance two competing objectives: (Objective 1) ensuring high fidelity by minimizing the pointwise discrepancy between the true solution and data-twin model prediction, and (Objective 2) preserving structural alignment by maximizing the normalized cosine similarity or correlation between their respective dynamics. 

The optimization problem (\ref{optimprob}) is resolved numerically using Pareto front analysis \citep{deb2001multiobjective}, a multi-objective optimization strategy that identifies the set of non-dominated solutions balancing competing objectives. The implementation details are provided in Appendix~\ref{appB}.

\subsection{Online Phase: Modeling the Data-Twin Temporal Dynamic by Explainable Deep Learning}\label{onphase}

The calculation step of the modes in the offline phase is performed
only once. The data-twin model is not complete without the temporal component.
These computations are carried out in the second step, during which the mathematical model is integrated over time. Additionally, if needed, the dynamics of the reduced-order model can be evaluated at any specific time point. For this reason, this step is referred to as the online phase, as it can be executed on demand to simulate the behavior of the original process.

Unlike conventional machine learning models, which typically function as opaque black boxes operating on entire datasets without disentangling underlying factors, a key innovation of this study lies in the explicit modeling of temporal dynamics decoupled from spatial dependencies. This is achieved through the use of explainable deep learning techniques, which offer interpretable insights into the data-twin model's parameters.  

In this work, Nonlinear AutoRegressive models with eXogenous inputs (NLARX) are employed as a novel deep learning-based approach within the domain of nonlinear system identification.
Since the emergence of artificial neural networks as numerical tools \citep{Judi1995, Ljung1999}, NLARX models have been used for various purposes, ranging from simulation \cite{Nelles2001}, to nonlinear predictive control \cite{Liu1998} or higher order nonlinear optimization problems \cite{Tieleman2012, Wang2018}.

The objective is to identify a dynamical system that effectively captures the nonlinear temporal evolution of the underlying process, expressed in the form:
\begin{equation}\label{tempdyn}
\left\{ \begin{array}{l}
\frac{{d{{\hat a}_j}\left( t \right)}}{{dt}} = {\mathcal{M}}\left( {{{\hat a}_j}\left( t \right),{a_j}\left( t \right),t} \right),\;\quad t \in {\mathbb{R}_{ \ge 0}}\\
{{\hat a}_j}\left( {{t_0}} \right) = \hat a_j^0,
\end{array} \right.
\end{equation}
based on modal temporal coefficients ${a_j}\left( {{t_i}} \right),\;j = 1,...,{N_{DTM}},\;\;{t_i} \in \left\{ {{{\rm{t}}_1},...,{{\rm{t}}_{{N_t}}}} \right\},$ computed by KROD algorithm in the offline phase, and initial condition obtained by projecting the first snapshot onto the Koopman mode basis:
\begin{equation}\label{inicond}
\hat a_j^0 = \left\langle {{u_0},{\phi _j}} \right\rangle ,j = 1,...,{N_{DTM}}.
\end{equation}

To approximate the system dynamics, the evolution of the state $\hat a\left( t \right)$ is modeled using the NLARX model structure, which is expressed as:
\begin{equation}\label{nlarx}
\hat a\left( t \right) = \mathcal{F}\left[ {\hat a\left( {t - 1} \right),...,\hat a\left( {t - {n_a}} \right),a\left( {t - {n_k}} \right),...,a\left( {t - {n_k} - {n_b} + 1} \right)} \right] + e\left( t \right),
\end{equation}
where $\mathcal{F}$ is a nonlinear estimator realized by a cascade forward neural network, capturing the underlying nonlinear temporal dependencies and input effects.
Here, ${{n_a}}$ denotes the number of past output terms, ${{n_b}}$ represents the number of past input terms used to predict the current output, and ${{n_k}}$ corresponds to the input delay, while $e\left( t \right)$ captures the modeling error. Each output of the NLARX model described in (\ref{nlarx}) is determined by regressors, which are transformations of previous input and output values. Typically, this mapping consists of a linear component and a nonlinear component, with the final model output being the sum of these two contributions.

Thus, the NLARX model serves as an empirical approximation of the continuous dynamics $\mathcal{M}$, with the neural network mapping $\mathcal{F}$ implicitly representing the evolution operator of the system, such that:
\begin{equation}\label{disnlarx}
\hat a\left( t \right) \approx \hat a\left( {t - 1} \right) + \Delta t \cdot \mathcal{M}\left( {\hat a\left( {t - 1} \right),a\left( {t - 1} \right),t - 1} \right) \approx \mathcal{F}\left[ {\hat a\left( {t - 1} \right),...,a\left(  \cdot  \right)} \right] + e\left( t \right),
\end{equation}
where $\Delta t$ is the discrete time step and $e\left( t \right)$ denotes the modeling error.

Training the NLARX model can be formulated as a nonlinear unconstrained optimization problem, of form:
\begin{equation}\label{nlarxoptim}
\theta \left( {{n_a},{n_b},{n_k}} \right) = \arg \min \,\frac{1}{{2{N_t}}}\sum\limits_{i = 1}^{{N_t}} {{{\left\| {a\left( {{t_i}} \right) - \widehat a\left( {{t_i}} \right)} \right\|}_2}}.
\end{equation}

In this formulation, the training dataset comprises the measured input $a\left( t \right)$, while $\widehat a\left( t \right)$ denotes the output generated by the NLARX model. The symbol ${\left\| {\, \cdot \,}\right\|_2}$ refers to the ${L_2}$ norm and $\theta \left( {{n_a},{n_b},{n_k}} \right)$ represents the parameter vector associated with the nonlinear function $\mathcal{F}$.

The NLARX architecture is well-suited for capturing system dynamics by feeding previous network outputs back into the input layer. It also allows for the specification of how many past input and output time steps are necessary to accurately represent the system's behavior. A critical aspect of effectively applying an NLARX network lies in the careful selection of inputs, input delays, and output delays. The objective is to adjust the network parameters $\theta \left( {{n_a},{n_b},{n_k}} \right)$ across the entire trajectory to minimize the objective function defined in (\ref{nlarxoptim}).

The reduced-order data-driven twin model for the PDE solution is obtained in the form:
\begin{equation}\label{phase2dtm}
\left\{ \begin{array}{l}
{u^{DTM}}\left( {x,t} \right) = \sum\limits_{j = 1}^{{N_{DTM}}} {{{\hat a}_j}\left( {{t}} \right){\phi _j}\left( x \right),} \;\;t \in {R_{ \ge 0}},\\
\hat a\left( t \right) = \mathcal{F}\left[ {\hat a\left( {t - 1} \right),...,\hat a\left( {t - {n_a}} \right),a\left( {t - {n_k}} \right),...,a\left( {t - {n_k} - {n_b} + 1} \right)} \right],\quad {{\hat a}_j}\left( {{t_0}} \right) = a_j^0,
\end{array} \right.
\end{equation}
where $N_{DTM}$ is the rank of the optimal Koopman basis, computed by the algorithm, ${\phi _j}\left( x \right),j = 1,...,{N_{DTM}}$ are the leading Koopman modes representing dominant spatial patterns computed by the algorithm, $\hat a\left( t \right)$ is the system temporal dynamics in the reduced-order space, modeled by the deep learning architecture with initial condition given by Eq.(\ref{inicond}).

\subsection{Qualitative Analysis of the Data-driven Twin Model}\label{qanaly}

A qualitative assessment of the orthogonality between the Koopman modes derived using the KROD algorithm proposed in this study can be carried out through the Modal Assurance Criterion (MAC) \citep{Ewins2000}. 

The discrete MAC value between two mode vectors \( \phi_i, \phi_j \in \mathbb{C}^n \) is defined as:
\begin{equation}\label{ortval}
\mathrm{MAC}_{ij} = \frac{\left| \left\langle \phi_i, \phi_j \right\rangle \right|^2}{\left\langle \phi_i, \phi_i \right\rangle \cdot \left\langle \phi_j, \phi_j \right\rangle},
\end{equation}
where \( \left\langle \phi_i, \phi_j \right\rangle = \phi_i^H \phi_j \) denotes the Hermitian (complex-conjugate) inner product in \( \mathbb{C}^n \). 

The orthogonality matrix \( \mathrm{MAC} \in \mathbb{R}^{N_{DTM} \times N_{DTM}} \) is symmetric and real-valued. The $MA{C_{ij}}$ values lie within the interval $\left[ {0,\;1} \right]$ as the leading Koopman modes obtained from the KROD algorithm are normalized vectors. In this context, Koopman modes are considered orthogonal if:
\begin{equation}\label{ort}
MA{C_{ij}} = \left\{ \begin{array}{l}
1,\quad i = j\\
0,\quad i \ne j.
\end{array} \right.
\end{equation}

The Pearson correlation coefficient is employed to quantify the correlation between the data-twin model solution and the true solution of the partial differential equation, providing a measure of how closely the data-twin model approximates the true dynamics across the space-time domain.

Let \( u(x, t) \) be the true solution and \( {u^{DTM}}\left( {x,t} \right) \) be the data-twin model solution of PDE, sampled at \( N_x \) spatial points and \( N_t \) time instances. Define the flattened vectors:

\[
\mathbf{U} = {[u({x_1},{t_1}), \ldots ,u({x_{{N_x}}},{t_1}),u({x_1},{t_2}), \ldots ,u({x_M},{t_{{N_t}}})]} \in {\mathbb{R}^{{N_x} {N_t}}},
\]
\[
\mathbf{U^{DTM}} = {[{u^{DTM}}({x_1},{t_1}), \ldots ,{u^{DTM}}({x_{{N_x}}},{t_1}),{u^{DTM}}({x_1},{t_2}), \ldots ,{u^{DTM}}({x_M},{t_{{N_t}}})]} \in {\mathbb{R}^{{N_x} {N_t}}}.
\]

Let \( U_k \) denote the \( k\text{-th} \) element of the vector \( \mathbf{U} \), and \( U_k^{\mathrm{DTM}} \) denote the \( k\text{-th} \) element of the vector \( \mathbf{U^{DTM}} \), respectively.

The Pearson correlation coefficient \( \rho \) between \( \mathbf{U} \) and \( \mathbf{U_{DTM}} \) is given by:
\begin{equation}\label{corr}
\rho  = \frac{{\sum\limits_{k = 1}^{{N_x}{N_t}} {\left( {{U_k} - \bar U} \right)} \left( {U_k^{DTM} - \overline {{U^{DTM}}} } \right)}}{{\sqrt {\sum\limits_{k = 1}^{{N_x}{N_t}} {{{\left( {{U_k} - \bar U} \right)}^2}} } \sqrt {\sum\limits_{k = 1}^{{N_x}{N_t}} {{{\left( {U_k^{DTM} - \overline {{U^{DTM}}} } \right)}^2}} } }},
\end{equation}
where the quantities
\[
\bar U = \frac{1}{{{N_x}{N_t}}}\sum\limits_{k = 1}^{{N_x}{N_t}} {{U_k}} ,\quad \overline {{U^{DTM}}}  = \frac{1}{{{N_x}{N_t}}}\sum\limits_{k = 1}^{{N_x}{N_t}} {U_k^{DTM}} 
\]
represent the mean values of the vectors \( \mathbf{U} \) and \( \mathbf{U}^{\mathrm{DTM}} \), respectively.

In addition to the Pearson correlation coefficient, the mean absolute error (MAE) is employed as a complementary metric to assess the average magnitude of pointwise discrepancies over the space-time domain, providing a more direct measure of predictive accuracy irrespective of correlation structure:
\begin{equation}\label{eror}
MAE = \frac{1}{{{N_x}{N_t}}}\sum\limits_{k = 1}^{{N_x}{N_t}} {\left| {{U_k} - U_k^{DTM}} \right|}.
\end{equation}

In addition to these metrics, to provide fine-grained, pointwise insight into how well the data-twin model approximates the true solution of PDE at specific locations in space and time, we define the absolute local error matrix $E \in {\mathbb{R}^{{N_x} \times {N_t}}}$ as:
\begin{equation}\label{errormat}
{E_{ij}} = \left| {u\left( {{x_i},{t_j}} \right) - {u^{DTM}}\left( {{x_i},{t_j}} \right)} \right|
\end{equation}
for all $i = 1,...,{N_x}$ and $j = 1,...,{N_t}$.

\subsection{Time Simulation and Validation of the Data-driven Twin Model}\label{timsim}

Data-driven twin model is subjected to time-domain simulation, followed by a two-step input-output validation of its predicted response. In the first validation step, the model is trained using the initial two-thirds of the available snapshots, with the remaining one-third reserved for validation. In the subsequent step, the training dataset is further reduced to the first one-third of the snapshots, and the remaining two-thirds are employed to assess the model's predictive performance.

\section{Data-driven Twin Modeling of Shock Wave Phenomena Using KROD}

This section demonstrates the application of the proposed KROD algorithm to shock wave modeling using the viscous Burgers' equation.

\subsection{Governing Equations of the Mathematical Model}\label{model}

The viscous Burgers equation model is considered, of the form:
\begin{equation}\label{burgers}
\left\{ \begin{array}{l}
\frac{\partial }{{\partial t}}u\left( {x,t} \right) + \frac{\partial }{{\partial x}}\left( {\frac{{u{{\left( {x,t} \right)}^2}}}{2}} \right) = \nu \frac{{{\partial ^2}}}{{\partial {x^2}}}u\left( {x,t} \right),\quad t > 0,\quad \nu  > 0,\\
u\left( {x,0} \right) = {u_0}\left( x \right),\quad x \in \mathbb{R},
\end{array} \right.
\end{equation}
where $u\left( {x,t} \right)$ is the unknown function of time $t$, $\nu $ is the viscosity parameter.

Three experiments are considered, exhibiting progressively more complex dynamics, each characterized by an initial condition of the following form:
\begin{equation}\label{test1}
Experiment\;1:\quad {u_0}\left( x \right) =  - \sin \left( {\pi x} \right),
\end{equation}
\begin{equation}\label{test2}
Experiment\;2:\quad {u_0}\left( x \right) = \left\{ \begin{array}{l}
{u_L},\quad x \le 0\\
{u_R},\quad x > 0
\end{array} \right.,
\end{equation}
\begin{equation}\label{test3}
Experiment\;3:\quad {u_0}\left( x \right) =  - \cos {(1.5\pi x)^2}.
\end{equation}

Homogeneous Dirichlet boundary conditions are imposed for all three examples, specified as follows:
\begin{equation}\label{boundcond}
u\left( {0,t} \right) = u\left( {L,t} \right) = 0
\end{equation}

Equation~(\ref{test1}), which defines the initial condition for Experiment 1, produces a sinusoidal pulse characterized by an abrupt change in slope at the boundaries of the domain.  
The initial value problem with the discontinuous initial condition given by (\ref{test2}) is referred to as a Riemann problem, which leads to the formation of a shock wave.   
For Experiment 3, the initial condition described by (\ref{test3}) induces a more complex evolution of the pulse.  

All three experiments involve phenomena that pose significant challenges for standard numerical methods to accurately capture.

\subsection{Analytical and Numerical Derivation of the Exact Solution}\label{annumsol}

The Cole-Hopf transformation, independently introduced by Hopf (1950) \citep{Hopf1950} and Cole (1951) \citep{Cole1951}, converts the Burgers' equation into the heat equation, thereby enabling the derivation of its exact solution. 
The analytical solution of the Burgers equation model (\ref{burgers}) is derived using the Cole–Hopf transformation, as detailed in Appendix \ref{colehopf}.

Using the Cole-Hopf transformation, the analytic solution to the problem (\ref{burgers}) is obtained in the following form:
\begin{equation}\label{exactsolu}
u\left( {x,t} \right) = \frac{{\int_{ - \infty }^\infty  {\frac{{x - \xi }}{t}{\varphi _0}\left( \xi  \right){e^{ - \frac{{{{\left( {x - \xi } \right)}^2}}}{{4\nu t}}}}d\xi } }}{{\int_{ - \infty }^\infty  {{\varphi _0}\left( \xi  \right){e^{ - \frac{{{{\left( {x - \xi } \right)}^2}}}{{4\nu t}}}}d\xi } }}.
\end{equation}

To numerically compute the exact solution (\ref{exactsolu}) corresponding to the initial conditions of the three test cases (\ref{test1})–(\ref{test3}), the Gauss–Hermite quadrature technique is employed \citep{Brass2011}, which is fully described in Appendix \ref{gauss}. 
The following results are thereby obtained for the three experimental cases.

The exact solution to the Burgers equation model (\ref{burgers}) with the initial condition (\ref{test1}) is equivalent with the following form:
\begin{equation}\label{exsol1}
u\left( {x,t} \right) = \frac{{\int_{ - \infty }^\infty  {4\nu z\;{e^{ - \frac{1}{{2\nu \pi }}\cos \left[ {\pi \left( {x - z\sqrt {4\nu t} } \right)} \right]}}{e^{ - {z^2}}}dz} }}{{\int_{ - \infty }^\infty  {\sqrt {4\nu t} \;{e^{ - \frac{1}{{2\nu \pi }}\cos \left[ {\pi \left( {x - z\sqrt {4\nu t} } \right)} \right]}}{e^{ - {z^2}}}dz} }}.
\end{equation}

In the case of Experiment 2 the exact solution to the Burgers equation model (\ref{burgers}) with the initial condition (\ref{test2}) is:
\begin{equation}\label{exsol2}
u\left( {x,t} \right) = \frac{{\int_{ - \infty }^\infty  {4\nu z\;{e^{ - \frac{{{u_R}}}{{2\nu }}\left( {x - z\sqrt {4\nu t} } \right)}}{e^{ - {z^2}}}dz} }}{{\int_{ - \infty }^\infty  {\sqrt {4\nu t} \;{e^{ - \frac{{{u_R}}}{{2\nu }}\left( {x - z\sqrt {4\nu t} } \right)}}{e^{ - {z^2}}}dz} }}.
\end{equation}

In the case of Experiment 3 the exact solution to the Burgers equation model (\ref{burgers}) with the initial condition (\ref{test3}) is:
\begin{equation}\label{exsol3}
u\left( {x,t} \right) = \frac{{\int_{ - \infty }^\infty  {4\nu z\;{e^{\frac{1}{{4\nu }}\left( {\frac{1}{{3\pi }}\sin \left( {3\pi (x - z\sqrt {4\nu t} )} \right) + x - z\sqrt {4\nu t} } \right)}}{e^{ - {z^2}}}dz} }}{{\int_{ - \infty }^\infty  {\sqrt {4\nu t} \;{e^{\frac{1}{{4\nu }}\left( {\frac{1}{{3\pi }}\sin \left( {3\pi (x - z\sqrt {4\nu t} )} \right) + x - z\sqrt {4\nu t} } \right)}}{e^{ - {z^2}}}dz} }}.
\end{equation}

\section{Numerical Results}\label{numres}

The algorithm underlying the Koopman Randomized Orthogonal Decomposition method is described in detail in Appendix \ref{appB}. This section presents numerical results highlighting the computational efficiency of the proposed KROD method, combined with Deep Learning (DL) for temporal simulation. The findings demonstrate that the approach produces a high-fidelity, reduced-complexity model that accurately captures the dynamics of the original system and enables real-time simulation, effectively functioning as a data-driven twin model.

The computational domain is defined as $\left[ {0,L} \right]$, with $L=2$, and the time interval considered is $\left[ {0,T} \right]$, where $T=3$. The viscosity parameter in the Burgers equation is set to $\nu  = {10^{ - 2}}$ for all three experiments. The domain is uniformly discretized using $N=100$ grid points, resulting in a mesh size of $\Delta x = 0.02$.
  
The Riemann problem parameters are assumed to be \( u_L = 0.1 \) and \( u_R = 0.5 \).  

By applying the Cole-Hopf transformation (\ref{cole}), the analytical solutions to the Burgers problem (\ref{burgers}) are obtained for all three experiments, corresponding respectively to the forms (\ref{exsol1}), (\ref{exsol2}), and (\ref{exsol3}), as detailed in Section \ref{cole}.  
Subsequently, the exact solutions for the three test cases are computed numerically using the Gauss-Hermite quadrature method (\ref{gh}) with \( n = 100 \) nodes (see Section \ref{gauss}).  

\begin{figure}[H]
\centerline{\includegraphics[width=12 cm]{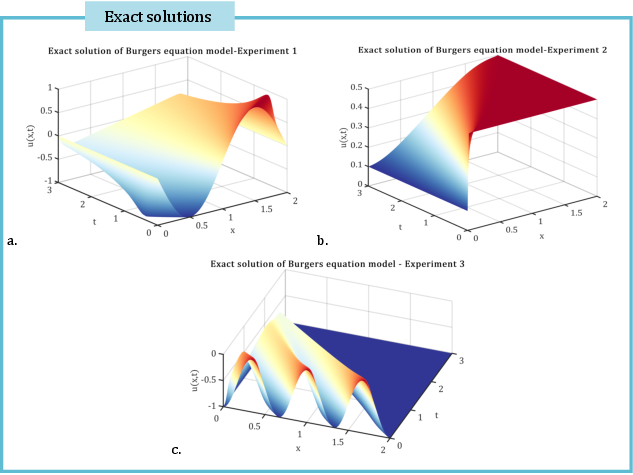}}
\caption{a) Exact solution of the viscous Burgers equation for the initial condition in Eq.~(\ref{test1}), derived analytically as in Eq.~(\ref{exsol1}); b) Exact solution of the viscous Burgers' equation for the second experiment with initial condition from Eq.~(\ref{test2}), derived analytically as in Eq.~(\ref{exsol2}); c) Exact solution of the viscous Burgers equation for the third experiment with initial condition from Eq.~(\ref{test3}), derived analytically as in Eq.~(\ref{exsol3}).\label{exact}}
\end{figure}

The following figure presents the exact solution to the viscous Burgers equation model.
Figure~\ref{exact}a presents the exact solution of the viscous Burgers' equation corresponding to the initial condition given in Eq.~(\ref{test1}), as derived in the analytical form shown in Eq.~(\ref{exsol1}). The exact solution of the viscous Burgers' equation for the second experiment, based on the initial condition specified in Eq.~(\ref{test2}), is depicted in Figure~\ref{exact}b. This solution is obtained in the analytical form provided in Eq.~(\ref{exsol2}). Figure~\ref{exact}c shows the exact solution of the viscous Burgers' equation for the third experiment, corresponding to the initial condition defined in Eq.~(\ref{test3}), and derived in the analytical form given in Eq.~(\ref{exsol3}).

The data-driven twin models of the form (\ref{phase1dtm}) are generated using the KROD algorithm as high-fidelity surrogates for the three Burgers problems examined in this study.  
For this purpose, the training dataset consists of a total of \( N_t  = 300 \) snapshots collected at uniformly spaced time intervals of \(\Delta t = 0.01\),  
and \( N_x = 101 \) spatial measurements recorded for each time snapshot.

KROD algorithm determines the optimal dimension \(N_{DTM}\) of the subspace spanned by the dominant Koopman modes by solving the multiobjective  optimization problem with nonlinear constraints stated in (\ref{optimprob}).  
This problem is solved by employing a genetic algorithm to identify the Pareto front corresponding to the two fitness functions.  

\begin{figure}[H]
\centerline{\includegraphics[width=12 cm]{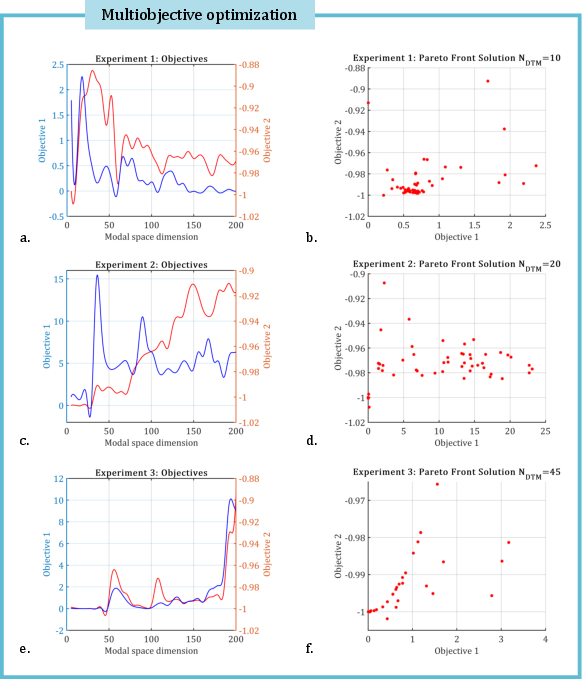}}
\caption{Figures (a), (c), and (e) show the optimization objectives defined in Eq.~(\ref{optimprob}) with respect to the modal space dimension in the range $\left[ {5,200} \right]$. Figures (b), (d), and (f) present the corresponding Pareto front solutions for the three experiments, obtained using a genetic algorithm.\label{Pareto}}
\end{figure}

Figures~\ref{Pareto} (a), (c) and (e) depict the optimization objectives defined in Eq.~(\ref{optimprob}), with respect to the modal space dimension in the range $\left[ {5,200} \right]$.
Figures~\ref{Pareto} (b), (d) and (f) present the Pareto front solution for each of the three experiments, respectively.

\begin{figure}[H]
\centerline{\includegraphics[width=10 cm]{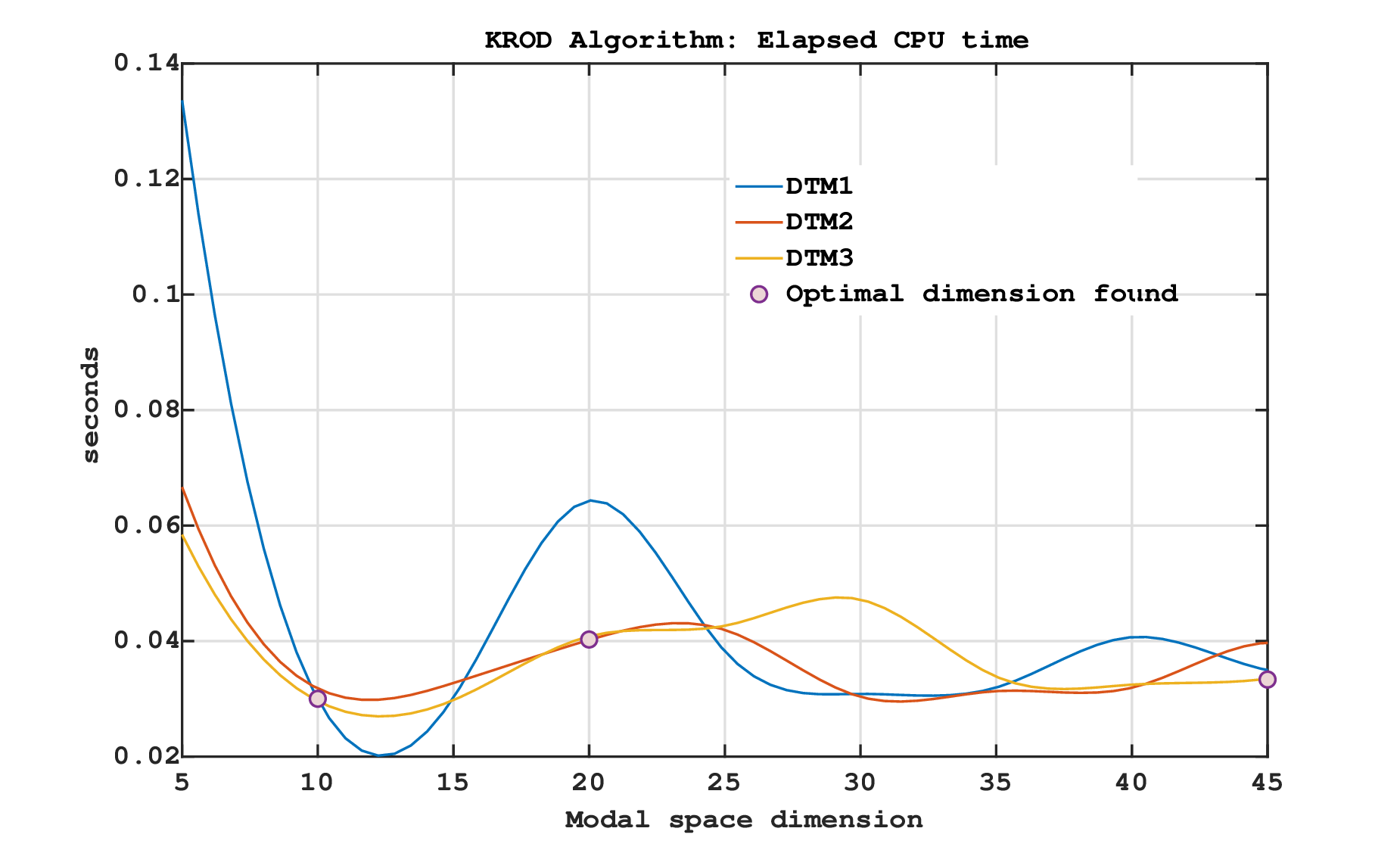}}
\caption{The running CPU time for KROD in offline phase, for the three experiments, respectively. \label{CPU}}
\end{figure}

For Experiment 1, the optimal dimension of the subspace spanned by the leading modes is \( N_{DTM} = 10 \) (Figure~\ref{Pareto}b).  
In Experiment 2, the data-twin model representing the Riemann problem dynamics has an optimal subspace dimension of \( N_{DTM} = 20 \) (Figure~\ref{Pareto}d).  
For Experiment 3, the increased complexity of the pulse is captured by a data-twin model with a higher dimension of \( N_{DTM} = 45 \) (Figure~\ref{Pareto}f).  
In all three cases, the computational time required by the KROD algorithm in offline phase remains extremely low, as illustrated in Figure~\ref{CPU}.

KROD algorithm identifies the optimal dimension of the modal basis and computes the optimal Koopman modes along with their associated time-dependent amplitudes.  
As a result, the data-driven twin model (\ref{phase1dtm}) is constructed, which-by definition-exhibits minimal approximation error and maximal correlation with respect to the original data model.

To facilitate a qualitative assessment of the DTMs, Table~\ref{qual1} reports the outcomes of the multiobjective optimization problem (\ref{optimprob}) as solved by the KROD.  
Specifically, the table includes the dimensionality of the leading Koopman modes, the absolute error as defined in Eq.~(\ref{eror}), and the correlation coefficient as given by Eq.~(\ref{corr}) between the exact solution and the corresponding data-twin model for each of the three experiments considered.

An analysis of the results presented in Table~\ref{qual1} reveals that KROD algorithm produces models exhibiting perfect correlation with the original data, i.e., \( \rho  = 1 \), across all three test cases.  
Furthermore, the absolute errors remain remarkably low, on the order of \( MAE < \mathcal{O}(10^{-7}) \).
\begin{table}[ht]
\centering
\caption{Qualitative analysis of the DTMs.}
\label{qual1}
\begin{tabular}{lccc}
\hline
\textbf{Test case} & \textbf{DTM complexity\textsuperscript{1}} & \textbf{DTM Correlation\textsuperscript{2}} & \textbf{DTM Absolute error\textsuperscript{3}} \\
\hline
Experiment 1 & $N_{\text{DTM}} = 10$ & $1.0000$ & $8.7264 \times 10^{-7}$ \\
Experiment 2 & $N_{\text{DTM}} = 20$ & $1.0000$ & $6.4552 \times 10^{-8}$ \\
Experiment 3 & $N_{\text{DTM}} = 45$ & $1.0000$ & $2.9684 \times 10^{-7}$ \\
\hline
\end{tabular}

\vspace{0.5em}
\footnotesize{
\textsuperscript{1} Number of retained Koopman modes. \\
\textsuperscript{2} DTM Correlation given by Eq.~(\ref{corr}). \\
\textsuperscript{3} DTM Absolute error given by Eq.~(\ref{eror}).
}
\end{table}


The selection of modes solely based on their energy ranking proves effective only under specific conditions \citep{Noack2011,Noack2016}.  
Some modes, despite contributing minimally to the overall energy, may exhibit rapid temporal growth at low amplitudes or represent strongly damped high-amplitude structures.  
Such modes can play a significant role in enhancing the fidelity of the data-twin model.  
A key advantage of the KROD algorithm over traditional DMD-based approaches lies in its ability to generate a substantially lower-dimensional subspace that retains the most dynamically significant Koopman modes.  
As a result, the KROD method avoids the need for an additional-often subjective-mode selection criterion.  
Moreover, the orthogonality of the KROD modes contributes to their qualitative superiority, enabling accurate modeling with a reduced number of modes.

\begin{figure}[H]
\centerline{\includegraphics[width=12 cm]{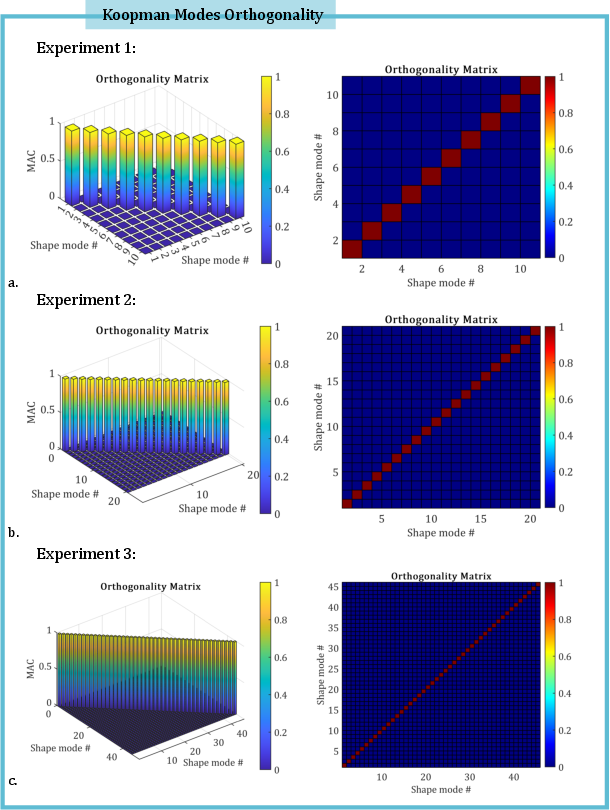}}
\caption{The orthogonality matrices corresponding to the three investigated experiments are shown, respectively, in both 3D and 2D views, in Figures (a), (b) and (c). The \( \mathrm{MAC} \) values, computed according to Eq.~(\ref{ortval}), confirm the mutual orthogonality of the Koopman modes in each case. \label{ortho}}
\end{figure}

To assess the orthogonality of the Koopman modes  computed by KROD algorithm, the orthogonality matrix introduced in Section~\ref{qanaly} is evaluated for each of the three experiments considered in this study.  
The orthogonality matrices corresponding to the three investigated experiments are shown, respectively, in both 3D and 2D views in Figure \ref{ortho}a–c.
The \( \mathrm{MAC} \) values calculated using Eq.~(\ref{ortval}) validate the mutual orthogonality of the Koopman modes across all cases.

\begin{figure}[H]
\centerline{\includegraphics[width=13 cm]{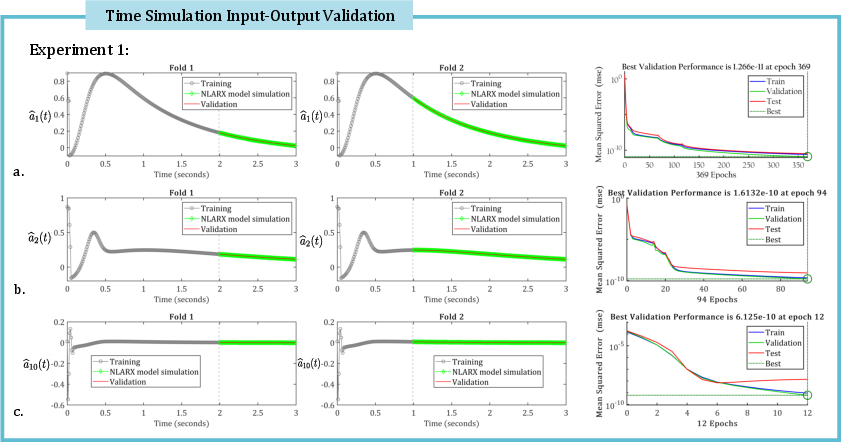}}
\caption{Two-fold input-output validation for the first two (a, b) temporal coefficients and the tenth (last) (c) temporal coefficient, during the simulated responses of the optimal NLARX models in Experiment~1. The corresponding deep learning performance metrics are also included in the figure on the right.\label{E1nlarx}}
\end{figure}

\begin{figure}[H]
\centerline{\includegraphics[width=13 cm]{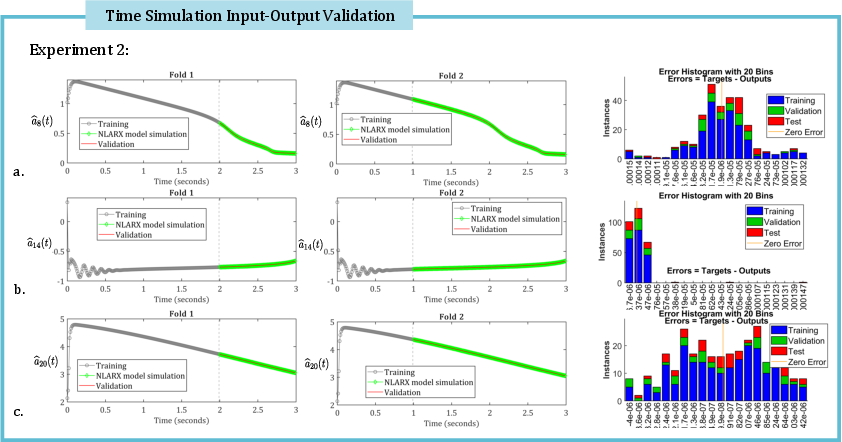}}
\caption{Two-fold input-output validation for the eighth (a) temporal coefficient, the fourteenth (b) temporal coefficient, and the twentieth (last) (c) temporal coefficient, during the simulated responses of the optimal NLARX models in Experiment~2. The corresponding deep learning performance metrics are also presented in the figure on the right. \label{E2nlarx}}
\end{figure}

\begin{figure}
\centerline{\includegraphics[width=13 cm]{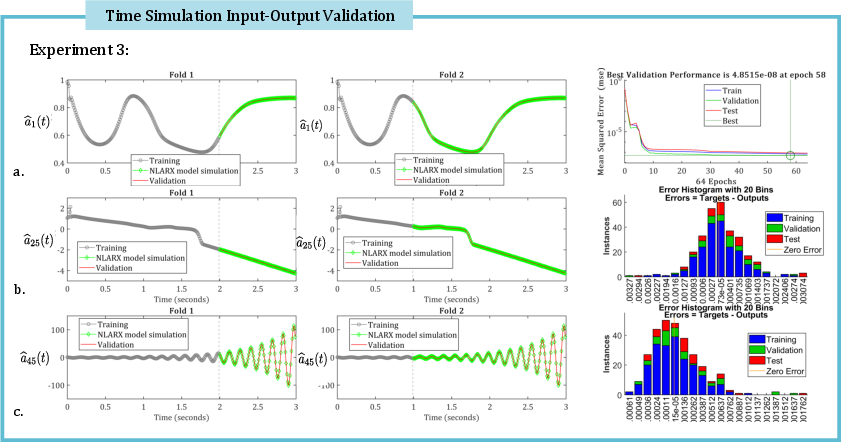}}
\caption{Two-fold input-output validation for the first (a) temporal coefficient, the twenty-fifth (b) temporal coefficient, and the forty-fifth (last) (c) temporal coefficient, during the simulated responses of the optimal NLARX models in Experiment~3. Corresponding deep learning performance metrics are presented in the figure on the right. \label{E3nlarx}}
\end{figure}

The KROD algorithm facilitates the identification of the leading Koopman modes and their corresponding temporal coefficients.
To achieve high-fidelity time simulation of the data-driven twin model's temporal parameters, Nonlinear AutoRegressive models with eXogenous inputs (NLARX) are employed (revise Subsection~(\ref{onphase}). 
The data-driven twin model, formulated as in Eq.~(\ref{phase2dtm}), is evaluated through time-domain simulation followed by a two-step input-output validation procedure. First, it is trained on the initial two-thirds of the snapshots, with the remaining third used for validation. Then, the training set is reduced to the first third, and the remaining two-thirds are used to evaluate predictive performance. 
Figures~\ref{E1nlarx}--\ref{E3nlarx} present the two-fold input-output validation results for three randomly selected temporal coefficients, simulated using the optimal NLARX models for Experiments 1, 2, and 3, respectively.  
These results demonstrate that the deep learning-based NLARX models accurately replicate the coefficient dynamics, even when the training dataset is reduced to only the first third of the available snapshots.  
Relevant performance metrics for the deep learning simulations are also provided in the figures.


\begin{figure}[H]
\centerline{\includegraphics[width=10 cm]{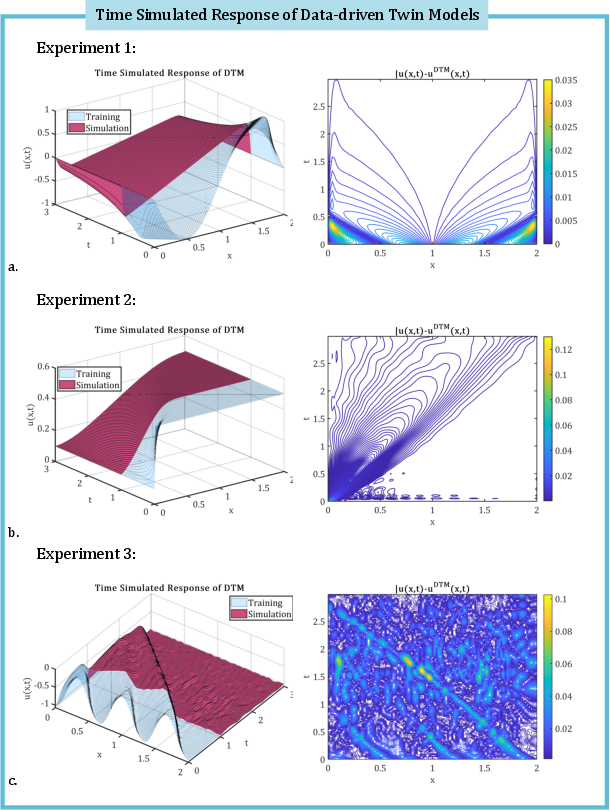}}
\caption{Data-driven twin model responses, constructed via Eq.~(\ref{phase2dtm}), show near-perfect agreement with exact solutions across all experiments, as evidenced by the low absolute local error from Eq.~(\ref{errormat}).
\label{simdtm}}
\end{figure}

Data-driven twin models, constructed in the form prescribed by Eq.~(\ref{phase2dtm}), are subsequently evaluated against the corresponding exact solutions to assess the absolute local error, as defined by Eq.~(\ref{errormat}), for each of the three experimental cases under consideration. As illustrated in Figure~\ref{simdtm}, the response of the twin models exhibits near-perfect correlation with the exact solutions, as evidenced by the consistently low magnitude of the local error.

\begin{figure}[H]
\centerline{\includegraphics[width=10 cm]{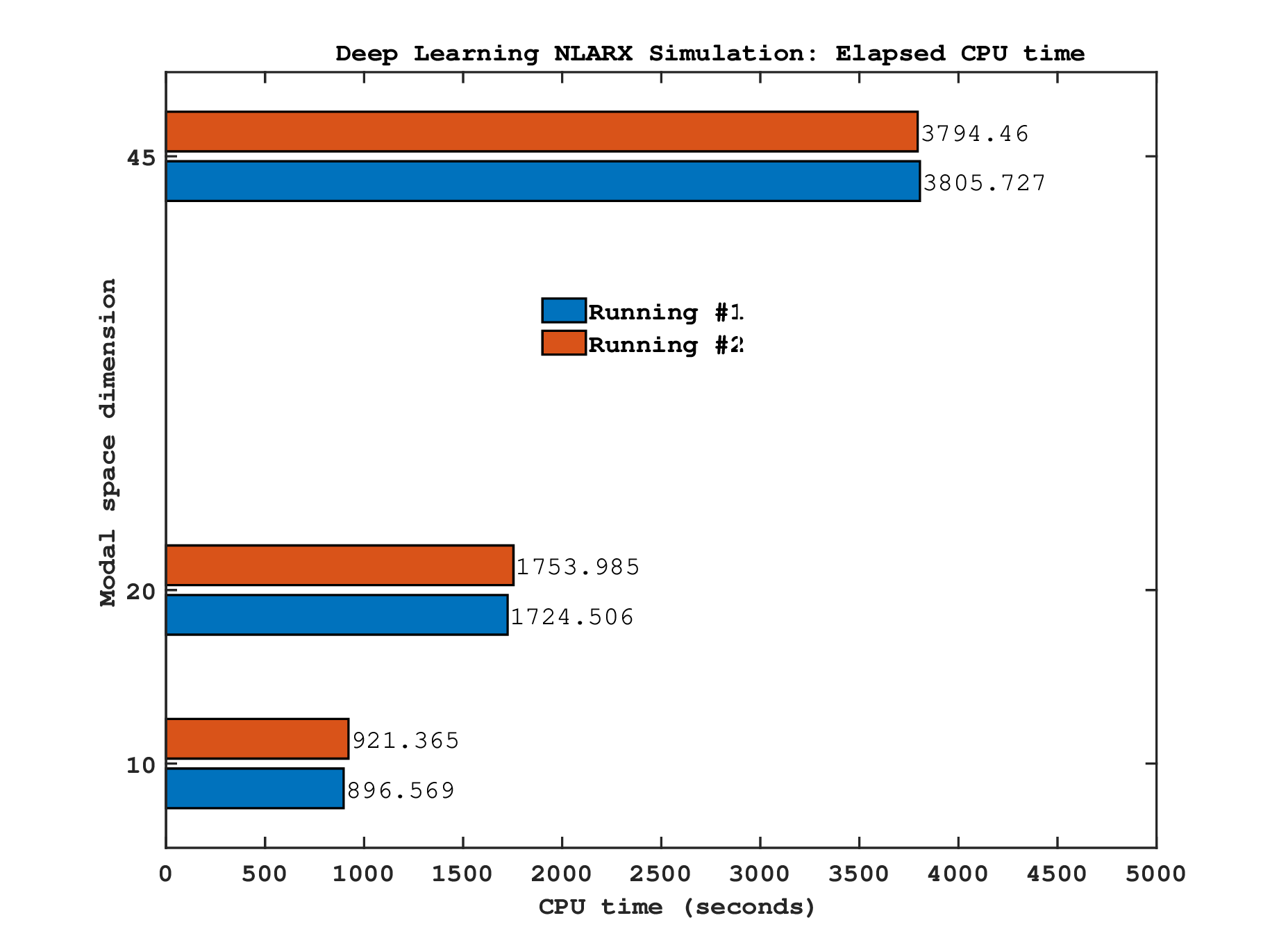}}
\caption{The running CPU time for the online phase of the algorithm, with respect to the modal space dimension optimized in the case of the three experiments considered.\label{DLCPU}}
\end{figure}

Figure \ref{DLCPU} presents the CPU time required during the online phase of the algorithm, in which the deep learning process is executed, as a function of the optimized modal space dimension for the three experimental cases analyzed. The reported CPU times correspond to two successive computational runs. All simulations and visualizations were performed using custom scripts developed in Matlab R2022a and executed on a system equipped with an Intel i7-7700K processor. The observed computational times are considered acceptable and fall within a reasonable range for practical implementation.

\section{Conclusions}

This study presents a novel framework for developing data-driven twin models with reduced computational complexity for nonlinear dynamical systems governed by partial differential equations, addressing a problem of growing importance in the field of data science. To the best of our knowledge, this is the first work to introduce the Koopman Randomized Orthogonal Decomposition (KROD) algorithm in the context of nonlinear PDEs.

This study introduces several methodological advancements that set it apart from classical modal decomposition techniques such as DMD and POD, which typically rely on heuristic or ad hoc criteria for mode selection. In contrast, the proposed KROD framework reduces computational cost while fully automating the identification of dominant modes.

The proposed approach comprises two primary phases. In the first phase, a novel numerical algorithm is developed to construct a reduced-order data-twin model, selected through Pareto front analysis to balance accuracy and complexity. In the second phase, an explainable NLARX deep learning framework is employed to enable real-time, adaptive calibration of the twin model and to predict its dynamic response. The resulting data-driven twin model achieves high fidelity by exhibiting minimal approximation error and strong correlation with the original dataset, thereby providing an accurate and computationally efficient surrogate for the underlying physical system.

A key novelty lies in the first quantitative framework for computing orthogonal Koopman modes using randomized orthogonal projections, which ensures optimal representation in the data space. Moreover, the proposed method integrates explainable deep learning techniques within the modal decomposition process, enabling interpretable and real-time calibration of the reduced-order data-twin model. 

Importantly, the proposed framework provides a pathway for future extensions into unsupervised learning paradigms, functioning effectively on raw observational data without requiring predefined models. This characteristic promotes broader applicability and generalization across a wide range of fluid dynamics problems and complex dynamical systems. The specific application of the KROD methodology to fluid dynamics scenarios will be the subject of future investigations.

In addition, the framework is readily extensible beyond the current scope. In particular, it holds promise for adaptation to stochastic partial differential equations, systems with sparse or incomplete measurements, and higher-dimensional problems involving 2D and 3D PDEs, directions the author aims to pursue in future research.

\vspace{6pt} 


\section*{Acknowledgments}
The author would like to express sincere gratitude to the anonymous reviewers for their in-depth and constructive feedback, which significantly contributed to improving the clarity and quality of this manuscript. The author also extends heartfelt thanks to Prof. Michael I. Navon for over a decade of invaluable collaboration and mentorship in the field of reduced-order modeling.

\section*{Abbreviations}

The following abbreviations are used in this manuscript:

\vspace{0.5em}

\begin{tabular}{ll}
\textbf{PDE}      & Partial Differential Equation \\
\textbf{DTM}      & Data-driven Twin Model \\
\textbf{POD}      & Proper Orthogonal Decomposition \\
\textbf{DMD}      & Dynamic Mode Decomposition \\
\textbf{RSVD}     & Randomized Singular Value Decomposition \\
\textbf{\textit{k}-RSVD} & Randomized Singular Value Decomposition of rank $k$ \\
\textbf{KROD}     & Koopman Randomized Orthogonal Decomposition \\
\textbf{NLARX}    & Nonlinear AutoRegressive models with eXogenous inputs \\
\textbf{MAC}      & Modal Assurance Criterion \\
\end{tabular}

\appendix

\section{Analytical Solution and Numerical Exact Solution }\label{appA}
\subsection{Derivation of the Analytical Solution Using the Cole-Hopf Transformation}\label{colehopf}

The Cole-Hopf transformation is defined by
\begin{equation}\label{cole}
u =  - 2\nu \frac{1}{\varphi }\frac{{\partial \varphi }}{{\partial x}}.
\end{equation}

Through analytical treatment, it is found that:
\begin{equation}\label{e1}
\frac{{\partial u}}{{\partial t}} = \frac{{2\nu }}{{{\varphi ^2}}}\left( {\frac{{\partial \varphi }}{{\partial t}}\frac{{\partial \varphi }}{{\partial x}} - \varphi \frac{{{\partial ^2}\varphi }}{{\partial x\partial t}}} \right),\quad u\frac{{\partial u}}{{\partial x}} = \frac{{4{\nu ^2}}}{{{\varphi ^3}}}\frac{{\partial \varphi }}{{\partial x}}\left( {\varphi \frac{{{\partial ^2}\varphi }}{{\partial {x^2}}} - \frac{{\partial \varphi }}{{\partial x}}\frac{{\partial \varphi }}{{\partial x}}} \right),
\end{equation}
\begin{equation}\label{e2}
\nu \frac{{{\partial ^2}u}}{{\partial {x^2}}} =  - \frac{{2{\nu ^2}}}{{{\varphi ^3}}}\left( {2{{\left( {\frac{{\partial \varphi }}{{\partial x}}} \right)}^3} - 3\varphi \frac{{{\partial ^2}\varphi }}{{\partial {x^2}}}\frac{{\partial \varphi }}{{\partial x}} + {\varphi ^2}\frac{{{\partial ^3}\varphi }}{{\partial {x^3}}}} \right).
\end{equation}

Substituting these expressions into (\ref{burgers}) it follows that:
\begin{equation}\label{e3}
\frac{{\partial \varphi }}{{\partial x}}\left( {\frac{{\partial \varphi }}{{\partial t}} - \nu \frac{{{\partial ^2}\varphi }}{{\partial {x^2}}}} \right) = \varphi \left( {\frac{{{\partial ^2}\varphi }}{{\partial x\partial t}} - \nu \frac{{{\partial ^3}\varphi }}{{\partial {x^3}}}} \right) = \varphi \frac{\partial }{{\partial x}}\left( {\frac{{\partial \varphi }}{{\partial t}} - \nu \frac{{{\partial ^2}\varphi }}{{\partial {x^2}}}} \right).
\end{equation}

Relation (A4) 
indicates that if $\varphi $ solves the heat equation, then $u\left( {x,t} \right)$ given by the Cole-Hopf transformation (A1) 
solves the viscous Burgers equation (\ref{burgers}). Thus  the viscous Burgers equation (\ref{burgers}) is recast into the following one:
\begin{equation}\label{heat}
\left\{ \begin{array}{l}
\frac{{\partial \varphi }}{{\partial t}} - \nu \frac{{{\partial ^2}\varphi }}{{\partial {x^2}}} = 0,\quad x \in R,t > 0,\nu  > 0,\\
\varphi \left( {x,0} \right) = {\varphi _0}\left( x \right) = {e^{ - \int_0^x {\frac{{{u_0}(\xi )}}{{2\nu }}d\xi } }},x \in \mathbb{R}.
\end{array} \right.
\end{equation}

Taking the Fourier transform with respect to $x$ for both heat equation and the initial condition (A5) 
it is found that:
\begin{equation}\label{solphi}
\varphi \left( {x,t} \right) = \frac{1}{{2\sqrt {\pi \nu t} }}\int\limits_{ - \infty }^\infty  {{\varphi _0}\left( \xi  \right)} \,{e^{ - \frac{{{{(x - \xi )}^2}}}{{4\nu t}}}}d\xi .
\end{equation}

Using the Cole-Hopf transformation (A1), 
the analytic solution to the problem (\ref{burgers}) is obtained in the following form:
\begin{equation}\label{exactsoluA}
u\left( {x,t} \right) = \frac{{\int_{ - \infty }^\infty  {\frac{{x - \xi }}{t}{\varphi _0}\left( \xi  \right){e^{ - \frac{{{{\left( {x - \xi } \right)}^2}}}{{4\nu t}}}}d\xi } }}{{\int_{ - \infty }^\infty  {{\varphi _0}\left( \xi  \right){e^{ - \frac{{{{\left( {x - \xi } \right)}^2}}}{{4\nu t}}}}d\xi } }}.
\end{equation}


\subsection{Derivation of the Numerical Exact Solution Using the Gauss-Hermite Quadrature}\label{gauss}

To compute the exact solution (A7) 
corresponding to the initial conditions of the three test cases (\ref{test1})–(\ref{test3}), the Gauss-Hermite quadrature technique \citep{Brass2011} is employed.

Gauss–Hermite quadrature approximates the value of integrals of the following kind:
\begin{equation}\label{gh}
\int\limits_{ - \infty }^\infty  {f\left( z \right)} \,{e^{ - {z^2}}}dz \approx \sum\limits_{i = 1}^n {{w_i}f\left( {{x_i}} \right)} ,
\end{equation}
where $n$ represents the number of sample points used, ${x_i}$ are the roots of the Hermite polynomial ${H_n}\left( x \right)$ and the associated weights ${w_i}$ are given by
\begin{equation}\label{weights}
{w_i} = \frac{{{2^{n - 1}}n!\sqrt \pi  }}{{{n^2}{{\left( {{H_{n - 1}}\left( {{x_i}} \right)} \right)}^2}}},\quad i = 1,...,n.
\end{equation}

In Experiment 1, the initial condition is defined  by Eq.(\ref{test1}):
\begin{equation}\label{ex1phi0}
{\varphi _0}\left( x \right) = {e^{ - \frac{1}{{2\nu }}\int_0^x {{u_0}(\xi )d\xi } }} = {e^{ - \frac{1}{{2\nu }}\int_0^x { - \sin \left( {\pi \xi } \right)d\xi } }} = {e^{\frac{1}{{2\nu \pi }}}} \cdot {e^{ - \frac{{\cos \left( {\pi x} \right)}}{{2\nu \pi }}}}.
\end{equation}
and the exact solution (A7) 
is obtained in the form:
\begin{equation}\label{sol1}
u\left( {x,t} \right) = \frac{{\int_{ - \infty }^\infty  {\frac{{x - \xi }}{t}{e^{ - \frac{{\cos \left( {\pi \xi } \right)}}{{2\nu \pi }}}} \cdot {e^{ - {{\left( {\frac{{x - \xi }}{{\sqrt {4\nu t} }}} \right)}^2}}}d\xi } }}{{\int_{ - \infty }^\infty  {{e^{ - \frac{{\cos \left( {\pi \xi } \right)}}{{2\nu \pi }}}} \cdot {e^{ - {{\left( {\frac{{x - \xi }}{{\sqrt {4\nu t} }}} \right)}^2}}}d\xi } }}.
\end{equation}

With the variable change:
\begin{equation}\label{varch}
z = \frac{{x - \xi }}{{\sqrt {4\nu t} }},
\end{equation}
the exact solution to the Burgers equation model (\ref{burgers}) with the initial condition (\ref{test1}) is equivalent with the following form:
\begin{equation}\label{exsol1A}
u\left( {x,t} \right) = \frac{{\int_{ - \infty }^\infty  {4\nu z\;{e^{ - \frac{1}{{2\nu \pi }}\cos \left[ {\pi \left( {x - z\sqrt {4\nu t} } \right)} \right]}}{e^{ - {z^2}}}dz} }}{{\int_{ - \infty }^\infty  {\sqrt {4\nu t} \;{e^{ - \frac{1}{{2\nu \pi }}\cos \left[ {\pi \left( {x - z\sqrt {4\nu t} } \right)} \right]}}{e^{ - {z^2}}}dz} }}.
\end{equation}

Similarly, by applying the variable transformation (A12), 
the exact solutions for the subsequent two experiments are derived as follows.

In the case of Experiment 2 the exact solution to the Burgers equation model (\ref{burgers}) with the initial condition (\ref{test2}) is:
\begin{equation}\label{exsol2A}
u\left( {x,t} \right) = \frac{{\int_{ - \infty }^\infty  {4\nu z\;{e^{ - \frac{{{u_R}}}{{2\nu }}\left( {x - z\sqrt {4\nu t} } \right)}}{e^{ - {z^2}}}dz} }}{{\int_{ - \infty }^\infty  {\sqrt {4\nu t} \;{e^{ - \frac{{{u_R}}}{{2\nu }}\left( {x - z\sqrt {4\nu t} } \right)}}{e^{ - {z^2}}}dz} }}.
\end{equation}

In the case of Experiment 3 the exact solution to the Burgers equation model (\ref{burgers}) with the initial condition (\ref{test3}) is:
\begin{equation}\label{exsol3A}
u\left( {x,t} \right) = \frac{{\int_{ - \infty }^\infty  {4\nu z\;{e^{\frac{1}{{4\nu }}\left( {\frac{1}{{3\pi }}\sin \left( {3\pi (x - z\sqrt {4\nu t} )} \right) + x - z\sqrt {4\nu t} } \right)}}{e^{ - {z^2}}}dz} }}{{\int_{ - \infty }^\infty  {\sqrt {4\nu t} \;{e^{\frac{1}{{4\nu }}\left( {\frac{1}{{3\pi }}\sin \left( {3\pi (x - z\sqrt {4\nu t} )} \right) + x - z\sqrt {4\nu t} } \right)}}{e^{ - {z^2}}}dz} }}.
\end{equation}

To obtain the numerical exact solutions for the three test cases, both the numerator and denominator are evaluated using Gauss-Hermite quadrature (A8). 

\section{Algorithms}\label{appB}

The following presents Algorithm 1, which outlines  the routine employed to compute the randomized singular value decomposition of rank $k$ ($k$-RSVD). Algorithm 2 details the Koopman Randomized Orthogonal Decomposition (KROD) procedure. 

\noindent\rule{15.8cm}{0.1pt}

\textbf{Algorithm 1: ($k$-RSVD) Randomized Singular Value Decomposition of Rank k }

\noindent\rule{15.8cm}{0.1pt}

\textbf{Initial data:} $V_0 \in {\mathbb{R}^{{N_x} \times {N_t}}}$, integer target rank $k \ge 2$ and $k \le {N_t}$.

\begin{enumerate}

\item[1.] Generate a Gaussian random test matrix $M \in \mathbb{R}^{N_t \times k}$ with independent and identically distributed (i.i.d.) standard normal entries.

\item[2.] Compute the compressed sampled range matrix $Q = V_0 M \in \mathbb{R}^{N_x \times k}$.

\item[3.] Produce Gram--Schmidt orthonormalisation of sampling matrix $Q$.

\item[4.] Form the projected matrix $P = Q^H V_0 \in \mathbb{R}^{k \times N_t}$.

\item[5.] Compute the economy-size singular value decomposition of $P$:
    \[
    P = T \Sigma W^H,
    \]
where $T \in \mathbb{R}^{k \times k}$ is orthogonal, $\Sigma \in \mathbb{R}^{k \times k}$ is diagonal, and $W \in \mathbb{R}^{N_t \times k}$ has orthonormal columns.

\item[6.] Compute the approximate left singular vectors as:
    \[
    U = Q T \in \mathbb{R}^{N_x \times k}.
    \]

\textbf{Output:}  Procedure returns $U \in {\mathbb{R}^{{N_x} \times k}}$, $\Sigma  \in {\mathbb{R}^{k \times k}}$, $W \in {\mathbb{R}^{{N_t}\times k}}$.

\end{enumerate}

\newpage
\noindent\rule{15.8cm}{0.1pt}

\textbf{Algorithm 2: (KROD) Koopman Randomized Orthogonal Decomposition}

\noindent\rule{15.8cm}{0.1pt}

\textbf{Initial data:} 
\begin{itemize}
\item A sequence of solution snapshots of PDE, arranged into two time-shifted data matrices:
\[{V_0} = \left[ {\begin{array}{*{20}{c}}
{{u_0}}&{{u_1}}&{...}&{{u_{{N_t} - 1}}}
\end{array}} \right]  \in {\mathbb{R}^{{N_x} \times {N_t}}},\;{V_1} = \left[ {\begin{array}{*{20}{c}}
{{u_1}}&{{u_2}}&{...}&{{u_{{N_t}}}}
\end{array}} \right]  \in {\mathbb{R}^{{N_x} \times {N_t}}}.\]

\item A sequence of integer target ranks: $\{k_1, k_2, \dots, k_N\}$ with $2 < k_1 < \dots < k_N \le N_t$.
\end{itemize}

\textbf{Offline phase:}
\begin{enumerate}
\item[1.] For $k \in \{k_1, k_2, \dots, k_N\}$:
  \begin{enumerate}
  \item Produce the randomized singular value approximation of rank $k$ for $V_0$: 
  \[V_0 \approx U \Sigma W^H.\]
The procedure for computing the $k$-RSVD is outlined in Algorithm 1.
  \item Compute the matrix operator $\mathcal{S}$ according to Eq.(\ref{defS}): 
  \[\mathcal{S} = {U^H}\left( {V_1W{\Sigma ^{ - 1}}} \right).\]
  \item Compute the Hermitian positive semi-definite Gram  matrix: 
  \[\mathcal{G} = {\mathcal{S}^H}\mathcal{S}.\]  
  \item Solve the eigenvalue problem $\mathcal{G}X=X\Lambda $ and obtain dynamic modes: $\Phi  = UX$. 
  \item Project \( V_0 \) onto the subspace spanned by the orthonormal basis \( \Phi \) to obtain the coefficient matrix: $A = \Phi^H V_0$.
  \end{enumerate}
  
\item[2.] Save corresponding model triplets: $\mathcal{T} = \{T_1, T_2, \dots, T_N\}$ where $T_i = (\Phi_i, A_i, k_i)$.

\item[3.] Multi-objective twin-model selection by solving the optimization problem (\ref{optimprob}) by Pareto front analysis: 
  \begin{enumerate}
  \item For each $T_i = (\Phi_i, A_i, k_i) \in \mathcal{T}$, 
    compute reduced-order approximation: 
    \[u^{(T_i)} = \sum_{j=1}^{k_i} a_j(t) \phi_j(x).\]
  \item Compute reconstruction error:
    \[E(T_i) = \| u - u^{(T_i)} \|_2.\]
  \item Compute normalized cosine similarity:
    \[C(T_i) = \frac{ \| u^H u^{(T_i)} \|_2^2 }{ \| u^H u \|_2 \cdot \| (u^{(T_i)})^H u^{(T_i)} \|_2 }.\]
  \item Construct the Pareto front $\mathcal{T}^*$ by selecting all non-dominated triplets:
\begin{align*}
\mathcal{T}^* = \big\{ T_i \in \mathcal{T} \;\big|\; \nexists T_j \in \mathcal{T} \text{ such that } &\; E(T_j) \le E(T_i), \\
& C(T_j) \ge C(T_i),\; \text{and strict in at least one} \big\}.
\end{align*}
  \item Return Pareto-optimal set of Koopman model triplet:
  \[\mathcal{T}^* = (\Phi_{\mathrm{DTM}}, A_{\mathrm{DTM}}, N_{\mathrm{DTM}}),\]
  consisting of the Koopman leading modes, the associated modal amplitudes and the rank of the Koopman basis, i.e. the number of terms in Koopman randomized orthogonal decomposition in form of Eq.(\ref{phase1dtm}).
  \end{enumerate}

\textbf{Online phase:}
\item[4.] Construct the nonlinear mapping for approximating the system dynamics by evolving the temporal coefficients via the optimal deep learning-based NLARX model Eq.~(\ref{disnlarx}).

\textbf{Output:} The reduced-order data-driven twin model for the PDE solution is obtained in the form given by Eq.(\ref{phase2dtm}).

\end{enumerate}

\end{document}